\documentclass[a4paper,draft,12pt,reqno]{article}
\usepackage{amsmath,amssymb,amsthm,xspace,amscd,latexsym,amssymb}
\usepackage[mathscr]{eucal}
\usepackage[all]{xy}
\usepackage{enumerate}

\theoremstyle{plain}
\newtheorem{theorem}{Theorem}
\newtheorem{example}[theorem]{Example}
\newtheorem{remark}[theorem]{Remark}
\newtheorem{lemma}[theorem]{Lemma}
\newtheorem{proposition}[theorem]{Proposition}
\newtheorem{corollary}[theorem]{Corollary}
\newtheorem{conjecture}[theorem]{Conjecture}

\newcommand{\tto}{\twoheadrightarrow}

\DeclareMathOperator{\Hom}{Hom}

\DeclareMathOperator{\Ext}{Ext}

\begin{document} 

\title{Some homological properties of the category $\mathcal{O}$}
\author{Volodymyr Mazorchuk}
\date{}
\maketitle

\begin{abstract}
In the first part of this paper the projective dimension of 
the structural modules in the BGG category $\mathcal{O}$ is 
studied. This dimension  is  computed for simple, standard and 
costandard modules.  For tilting and injective modules an explicit 
conjecture relating the result to Lusztig's $\mathbf{a}$-function
is formulated (and proved for type $A$). The second part deals
with the extension algebra of Verma modules. It is shown that
this algebra is in a natural way $\mathbb{Z}^2$-graded  and
that it has two $\mathbb{Z}$-graded Koszul subalgebras. The
dimension of the space $\mathrm{Ext}^1$ into the projective
Verma module is determined. In the last part several new classes
of Koszul modules and modules, represented by linear complexes
of tilting modules, are constructed.
\end{abstract}


\section{Introduction}\label{s0}

The Bernstein-Gelfand-Gelfand category $\mathcal{O}$, \cite{BGG}, 
associated with
a triangular decomposition of a semi-simple complex finite-dimensional
Lie algebra is an important and intensively studied object in 
modern representation theory. It has many very beautiful properties
and symmetries. For example it is equivalent to the module category
of a standard Koszul quasi-hereditary algebra and is Ringel self-dual.
Its principal block is even Koszul self-dual. 
Powerful tools for the study of the 
category  $\mathcal{O}$ are Kazhdan-Lusztig's combinatorics, developed 
in  \cite{KL}, and Soergel's combinatorics, worked out in \cite{So}. 
These two machineries immediately give a lot of information about the
numerical algebraic and homological invariants of simple, 
projective, Verma and tilting modules in $\mathcal{O}$ respectively. 
However, many natural questions about such invariants are
still open. The present paper answers some of them. 

The paper starts with a description of notation and preliminary
results in Section~\ref{s1}. The rest is divided into three parts.
The first part of this is Section~\ref{s2}, which is dedicated
to the study of homological dimension for structural modules
in the principal block $\mathcal{O}_0$ of $\mathcal{O}$. By {\em
structural} I mean projective, injective, simple, standard (Verma), 
costandard (dual Verma), and tilting modules respectively. 
In some cases the
result is rather expected. Some estimates go back to the original
paper \cite{BGG}. For simple and standard modules the result
can be deduced from Soergel's Koszul self-duality of $\mathcal{O}$.
However, to my big surprise I failed to find more elementary arguments
in the available literature. Here I present an explicit answer
for simple, standard and costandard modules, and a proof, which does
not even uses the Kazhdan-Lusztig conjecture. However,
the shortest ``elementary'' argument I could come up with uses some
properties of Arkhipov's twisting functors, established in \cite{AS}.
Things become really interesting when one tries to compute the
projective dimension of an indecomposable tilting module. Although the
projective dimension of the characteristic tilting module in 
$\mathcal{O}_0$ is well-known (see for example \cite{MO}), it seems
that nobody has tried to determine the projective dimension of an
indecomposable tilting module. A very surprizing conjecture based
on several examples and Theorem~\ref{t9}, which says that the
projective dimension of an indecomposable tilting module is a function,
constant on two-sided cells, suggests that this
dimension is given by Lusztig's $\mathbf{a}$-function from \cite{Lu}.
This conjecture is proved here for type $A$ (Theorem~\ref{t12}), which
might be considered as a good evidence that the result should be true
in general. However, I have no idea how to approach this question
in the general case and my arguments from type $A$ certainly
can't be transfered. The determination of the projective dimension
for injective modules reduces to that of tilting 
modules. As a ``bonus'' we also give a formula for the projective
dimension of Irving's shuffled Verma modules in Proposition~\ref{p15}.

In Section~\ref{s3} we study the extension algebra of
standard modules in $\mathcal{O}_0$. This is an old open problem,
where really not that much is known. The only available conjecture
about the numerical description of such extensions, formulated
in \cite[Section~5]{GJ}, is known to be false (\cite{Bo}), and the only
explicit partial results I was able to find is the ones obtained
in \cite{GJ,Ca}. Here I follow the philosophy of \cite{DM}, where
it was pointed out that the extension algebra of standard
modules is naturally $\mathbb{Z}^2$-graded. This $\mathbb{Z}^2$-grading 
is obtained from two different $\mathbb{Z}$-gradings: the first one which
comes from the category of graded modules, and the second one which
comes  from  the derived category. Koszul
self-duality of $\mathcal{O}_0$ induces a non-trivial automorphism
of this $\mathbb{Z}^2$-graded algebra, which swaps the
$\mathbb{Z}$-graded subalgebras of homomorphisms and linear
extensions, see Theorem~\ref{t2.02}. This allows one to
calculate linear extensions between standard modules, in
particular, to reprove the main result from \cite{Ca}. 
A surprizing corollary here is that by far not all 
projectives from the linear projective resolution of a
standard module give rise to a non-trivial linear extension
with the standard module, determined by this projective. 
In \cite{DM} it was shown that in the multiplicity
free case the extension algebra of standard modules is Koszul
(with respect to the $\mathbb{Z}$-grading, which is 
naturally induced by the $\mathbb{Z}^2$-grading mentioned above).
I do not think that this is true in the general case since I
do not believe that the extension algebra of standard modules
is generated in degree $1$. However, I think it is reasonable 
to expect that the subalgebra of this extension algebra,
generated by all elements of degree $1$, is Koszul. To support
this it is shown that the $\mathbb{Z}$-graded
subalgebra of all {\em homomorphisms} between standard modules is
Koszul, see Proposition~\ref{p2.05}. As the last result of 
Section~\ref{s3}  I explicitly determine
the dimension of the $\mathrm{Ext}^1$ space
from a standard module to a projective standard module,
see Theorem~\ref{t2.12}. From my point of view, the
answer is again surprizing.

In \cite{MO2,MOS} one finds an approach to Koszul duality
using the categories of linear complexes of projective
or tilting modules. For the category $\mathcal{O}_0$ this approach
can be used to get quite a lot of information, see
\cite{Ma,MO2,MOS}. In particular one can prove the Koszul
duality of various functors and various algebras, associated
to $\mathcal{O}_0$. A very important class of modules for
Koszul algebras is the class of the so-called {\em Koszul
modules}. These are modules with linear projective resolutions.
Such modules have a two-folded origin, namely, they
are both modules over the original algebra and over its Koszul
dual (via the corresponding linear resolution). In Section~\ref{s4} I 
show for several natural classes of modules from
$\mathcal{O}_0$ that they are either Koszul or can be
represented in the derived category by a linear complex
of tilting modules (which roughly means that they correspond
to Koszul modules for the Ringel dual of $\mathcal{O}_0$). 
The latter property seems to be more ``natural'' for the 
category $\mathcal{O}_0$. For example, while only the
simple and the standard modules are Koszul, it turns out that 
all simple, standard, costandard and shuffled Verma modules
are represented by linear complexes of tilting modules
(for the latter statement see Theorem~\ref{t16}). 
As an extension of this list we also show that some
structural modules from the parabolic subcategories
also have at least one of these properties, when considered as
objects in the original category $\mathcal{O}_0$.

\section{Notation and preliminaries}\label{s1}

Let $\mathfrak{g}$ denote a semi-simple finite-dimensional
Lie algebra over $\mathbb{C}$ with a fixed triangular decomposition,
$\mathfrak{g}=\mathfrak{n}_-\oplus\mathfrak{h}\oplus\mathfrak{n}_+$.
Let $\mathcal{O}$ denote the corresponding BGG-category $\mathcal{O}$,
defined in \cite{BGG}. Let $\mathcal{O}_0$ denote the {\em principal
block} of $\mathcal{O}$, that is the indecomposable direct summand of 
$\mathcal{O}$, containing the trivial module. Let $W$ be the Weyl group
of $\mathfrak{g}$ which acts on $\mathfrak{h}^*$ in the usual way 
$w(\lambda)$ and via the dot-action $w\cdot \lambda$. The category
$\mathcal{O}_0$ contains the Verma modules $M(w\cdot 0)$, $w\in W$.
For $w\in W$ we set $\Delta(w)=M(w\cdot 0)$ and let
$L(w)$ denote the unique simple quotient of $\Delta(w)$. Further,
$P(w)$ is the indecomposable projective cover of $L(w)$ and
$I(w)$ is the indecomposable injective envelope of $L(w)$. 
We set $L=\oplus_{w\in W}L(w)$ and analogously for all other
structural modules. 

The category $\mathcal{O}_0$ is a highest weight category in the
sense of \cite{CPS}, in particular, associated to $L(w)$
we also have the costandard module $\nabla(w)$, and the indecomposable 
tilting module $T(w)$ (see \cite{Ri}). If $\star$ is the standard 
duality on $\mathcal{O}$, we have $\nabla(w)\cong \Delta(w)^{\star}$
and $T(w)\cong T(w)^{\star}$.
For $w\in W$ by $l(w)$ we denote the length of $w$. 
Let $w_0$ denote the longest element of $W$. By $\leq$ we denote
the Bruhat order on $W$. For $w\in W$ let $\theta_w:\mathcal{O}_0\to\mathcal{O}_0$
denote the indecomposable projective functor, corresponding to
$w$, see \cite[Theorem~3.3]{BG}. In particular, if $s\in W$ is
a simple reflection, then $\theta_s$  is the
translation functor through the $s$-wall (see \cite[Section~3]{GJ}).
We have $\theta_wP(e)\cong P(w)$ (\cite[Theorem~3.3]{BG}) and
$\theta_wT(w_0)\cong T(w_0w)$ (\cite[Theorem~3.1]{CI}).

If $\mathcal{X}^{\bullet}$ is a complex and $n\in\mathbb{Z}$, by 
$\mathcal{X}^{\bullet}[n]$ we will denote the $n$-th shifted complex,
that is the complex, satisfying $(\mathcal{X}^{\bullet}[n])^{i}\cong
\mathcal{X}^{i+n}$ for all $i\in\mathbb{Z}$. We also use the standard
notation $\mathcal{D}^b(A)$, $\mathcal{L}\,\mathrm{F}$ and
$\mathcal{R}\,\mathrm{F}$ to denote the bounded derived category,
and the left and right derived functors respectively.

Let $A=\mathrm{End}_{\mathcal{O}}(P)^{\mathrm{op}}$ be the associative 
algebra of $\mathcal{O}_0$. This means that $\mathcal{O}_0$ is equivalent 
to the category $A\mathrm{-mod}$ of finitely generated left $A$-modules. 
This algebra is Koszul (\cite[Theorem~18]{So})
and we denote by $\mathtt{A}$ the associated positively graded
algebra. Denote by $\mathtt{A}\mathrm{-gmod}$ the category of all 
finitely generated  graded left $\mathtt{A}$-modules.
For $w\in W$ we denote by $\mathtt{L}(w)$ the standard 
graded lift of  $L(w)$, concentrated in degree $0$; and by 
$\mathtt{P}(w)$  and $\mathtt{I}(w)$ the corresponding lifts of  
$P(w)$ and $I(w)$ respectively  such that the maps $P(w)\tto L(w)$ 
and $L(w)\hookrightarrow I(w)$ become homogeneous of degree $0$.
Further we fix graded lifts $\mathtt{\Delta}(w)$ and 
$\mathtt{\nabla}(w)$ such that the obvious maps $P(w)\tto \Delta(w)$ 
and $\nabla(w)\hookrightarrow I(w)$ become homogeneous of 
degree $0$. Finally, we fix the graded lift $\mathtt{T}(w)$ such that
the map $\Delta(w)\hookrightarrow T(w)$ becomes homogeneous of 
degree $0$. In general, we will try to follow the conventions 
of  \cite[Introduction]{MOS} and refer the reader to that paper 
for details. In particular, a graded lift of a module, $M$, 
will be usually denoted by $\mathtt{M}$. For $k\in\mathbb{Z}$
we denote by $\langle k\rangle$ the functor of shifting
the grading as follows: if $\mathtt{M}=\oplus_{i\in\mathbb{Z}}
\mathtt{M}_i$ then $\mathtt{M}\langle k\rangle_i=\mathtt{M}_{i+k}$.
A complex $\mathcal{X}^{\bullet}$ of graded projective
(respectively injective or tilting) modules is called {\em linear} 
provided that $\mathcal{X}^{i}\in
\mathrm{add}(\mathtt{P}\langle i\rangle)$
(respectively $\mathtt{I}\langle i\rangle$ and
$\mathtt{T}\langle i\rangle$) for all $i\in\mathbb{Z}$. By 
$\mathscr{LC}(\mathtt{P})$ (respectively $\mathscr{LC}(\mathtt{I})$ 
or $\mathscr{LC}(\mathtt{T})$) we denote the category, whose objects
are all linear (bounded) complexes of projective (respectively 
injective and  tilting) modules, and morphisms are all possible 
morphisms of complexes of graded modules. For general information 
about the  categories of linear complexes and their applications, see 
\cite{MO2,MOS}.

\section{Projective dimensions of structural modules in 
$\mathcal{O}_0$}\label{s2}

As we already mentioned, the category $\mathcal{O}_0$ is a highest 
weight category. All simple, standard, costandard, projective, injective 
and tilting modules play various important roles in this structure. 
Our first natural question is to determine the projective dimension
of all these (indecomposable) structural modules. We will write
$\mathrm{p.d.}(M)$ for the {\em projective dimension} of a module, $M$,
and denote by $\mathrm{gl.dim.}$ the {\em global} (or
{\em homological}) {\em dimension} of an algebra or its module category.
As an
obvious result here one can mention $\mathrm{p.d.}(P(w))=0$
for all $w\in W$.

\subsection{Standard and simple modules}\label{s2.1}

It turns out that determining the projective dimension of 
standard and simple modules in $\mathcal{O}_0$ is the easiest
part of the task. Actually, first estimates for 
these dimensions were already obtained in the original paper 
\cite{BGG}. 

\begin{proposition}\label{p1}(\cite[Section~7]{BGG})
\begin{enumerate}[(i)]
\item\label{p1.1} $\mathrm{p.d.}(\Delta(w))\leq l(w)$.
\item\label{p1.2} $\mathrm{p.d.}(L(w))\leq 2l(w_0)-l(w)$.
\item\label{p1.3} $\mathrm{gl.dim.}\mathcal{O}_0\leq 2l(w_0)$.
\end{enumerate}
\end{proposition}

\begin{proof}
Obviously, $\mathrm{p.d.}(\Delta(e))=0$ since $\Delta(e)=P(e)$.
As we have already mentioned, $\mathcal{O}_0$ is a highest weight
category with respect to the Bruhat order on $W$. In particular, this
means that the kernel of the natural projection $P(w)\tto \Delta(w)$
has a filtration with subquotients $\Delta(w')$, $l(w')<l(w)$. Hence
\begin{displaymath}
\mathrm{p.d.}(\Delta(w))\leq 
\max_{w':l(w')<l(w)}\{\mathrm{p.d.}(\Delta(w'))\}+1,
\end{displaymath}
which implies \eqref{p1.1} by induction.

Since $\Delta(w_0)=L(w_0)$, the formula of \eqref{p1.2} for $w=w_0$ 
is just a special case of \eqref{p1.1}. Consider now the short exact sequence
$X\hookrightarrow \Delta(w)\tto L(w)$. Then $X$ has a filtration with
subquotients of the form $L(w')$, $l(w')>l(w)$. Hence one obtains 
\begin{displaymath}
\mathrm{p.d.}(L(w))\leq \max_{w':l(w')>l(w)}\{\mathrm{p.d.}(L(w'))\}+1,
\end{displaymath}
which implies \eqref{p1.2} by induction.

\eqref{p1.3} is an immediate corollary from \eqref{p1.2}.
\end{proof}

Further, in the last remark in \cite{BGG} it is mentioned that one can show
that $\mathrm{gl.dim.}\mathcal{O}_0= 2l(w_0)$. The shortest argument 
I know, which does this, is the following:

\begin{proposition}\label{p2}
$\mathrm{p.d.}(L(e))\geq 2l(w_0)$, in particular,
$\mathrm{gl.dim.}\mathcal{O}_0= 2l(w_0)$.
\end{proposition}

\begin{proof}
Consider the BGG-resolution
\begin{displaymath}
0\to M_{l(w_0)}\to M_{l(w_0)-1}\to\dots
\to M_{1}\to M_{0}\to L(e)\to 0
\end{displaymath}
of $L(e)$, see \cite[Theorem~10.1]{BGG2}, and let $\mathcal{M}^{\bullet}$ 
be the corresponding complex of (direct sums of) Verma modules, whose 
only non-zero homology is $\mathrm{H}^{0}(\mathcal{M}^{\bullet})\cong L(e)$. 
Every non-zero map $f:\Delta(w_0)\to\nabla(w_0)$ induces a non-zero
map $\overline{f}:\mathcal{M}^{\bullet}\to 
(\mathcal{M}^{\bullet})^{\star}[2l(w_0)]$. Since
$\dim\Hom_{\mathcal{O}}(\Delta(w),\nabla(w'))=\delta_{w,w'}$
by \cite[Section~3]{Ri}, it follows that $\overline{f}$ 
is not homotopic to $0$. Since
$\Ext_{\mathcal{O}}^i(\Delta(w),\nabla(w'))=0$ for all $i>0$
by \cite[Theorem~4]{Ri}, from \cite[Chapter~III(2), Lemma~2.1]{Ha}
it follows that
$\Ext_{\mathcal{O}}^{2l(w_0)}(L(e),L(e))\neq 0$. Thus
we get $\mathrm{p.d.}(L(e))\geq 2l(w_0)$. The latter and
Proposition~\ref{p1}\eqref{p1.3} imply
$\mathrm{gl.dim.}\mathcal{O}_0= 2l(w_0)$.
\end{proof}

Now let us show that the estimates in Proposition~\ref{p1}\eqref{p1.1}
and Proposition~\ref{p1}\eqref{p1.2} are in fact the exact values.
Already this becomes slightly tricky, especially for simple modules.
Here we present a uniform approach, which works for both standard
and simple modules, and is based on certain properties of the
so-called {\em twisting functors} on $\mathcal{O}_0$. 
Some other approaches will
be discussed in remarks at the end of this subsection. For
$w\in W$ let $\mathrm{T}_w:\mathcal{O}_0\to
\mathcal{O}_0$ denote the corresponding twisting functor,
see \cite{Ar,AS}. Let further $\mathrm{G}_{w^{-1}}:\mathcal{O}_0\to
\mathcal{O}_0$ denote the right adjoint of $\mathrm{T}_w$.
The functor $\mathrm{G}_{w^{-1}}$ is isomorphic 
(\cite[Corollary~6]{KM}) to Joseph's completion functor defined 
in \cite{Jo}. 
We start with  the case of standard modules since the proof is more 
direct in this case.

\begin{proposition}\label{p3}
$\Ext_{\mathcal{O}}^{l(w)}(\Delta(w),L(e))\neq 0$, in particular,
$\mathrm{p.d.}(\Delta(w))=l(w)$.
\end{proposition}

\begin{proof}
We do induction on $l(w)$. If $w=e$ the statement is
obvious. If $s$ is a simple reflection such that $l(sw)>l(s)$, we have

\begin{displaymath}
\begin{array}{lcl}
\mathrm{Ext}_{\mathcal{O}}^{l(sw)}(\Delta(sw),L(e)) & = & \\
\mathrm{Hom}_{\mathcal{D}^b(\mathcal{O})}(\Delta(sw),L(e)[l(sw)]) & = & 
\text{(by \cite[(2.3)]{AS})}\\
\mathrm{Hom}_{\mathcal{D}^b(\mathcal{O})}(\mathrm{T}_s\Delta(w),L(e)[l(sw)]) 
& = &  \text{(by \cite[Theorem~2.2]{AS})}\\
\mathrm{Hom}_{\mathcal{D}^b(\mathcal{O})}
(\mathcal{L}\mathrm{T}_s\Delta(w),L(e)[l(sw)]) 
& = &  \text{(by \cite[Corollary~4.2]{AS})}\\
\mathrm{Hom}_{\mathcal{D}^b(\mathcal{O})}
(\Delta(w),\mathcal{R}\mathrm{G}_sL(e)[l(sw)]) 
& = &  \text{(by \cite[Corollaries~4.2 and 6.2]{AS})}\\
\mathrm{Hom}_{\mathcal{D}^b(\mathcal{O})}
(\Delta(w),L(e)[l(sw)-1])& = &\\
\mathrm{Hom}_{\mathcal{D}^b(\mathcal{O})} (\Delta(w),L(e)[l(w)])& = &\\
\mathrm{Ext}_{\mathcal{O}}^{l(w)}(\Delta(w),L(e)) & \neq & 0
\end{array}
\end{displaymath}
by induction. The statement now follows from 
Proposition~\ref{p1}\eqref{p1.1}.
\end{proof}

\begin{remark}\label{r4}
{\em Another way to prove the formula for the projective dimension
of standard modules from 
Proposition~\ref{p3} is to use
\cite[Theorem~18]{So}, \cite[Proposition~2.7]{ADL} and 
\cite[3.5]{Ir}. A disadvantage in this case is the fact that 
so far there is no purely algebraic proof of \cite[Theorem~18]{So},
whereas the results from \cite{AS} used in the proof of
Proposition~\ref{p3} can be proved algebraically.
}
\end{remark}

\begin{remark}\label{r5}
{\em Yet another way to prove the formula for the projective dimension
of standard modules from  Proposition~\ref{p3} is to observe,
using translation functors, that  $\mathrm{p.d.}\Delta(w_0)$ 
coincides with the projective dimension of the characteristic 
tilting module in $\mathcal{O}_0$. Then \cite[Corollary~2]{MO}
and Proposition~\ref{p2} imply $\mathrm{p.d.}(\Delta(w_0))=l(w_0)$.
For any $w\in W$ and a simple reflection $s\in W$ such that
$l(ws)>l(w)$ there is a short exact sequence $\Delta(ws)\hookrightarrow
\theta_s\Delta(w)\tto\Delta(w)$. Since $\theta_s$ is exact
and maps projectives to projectives, we have
$\mathrm{p.d.}(\theta_s\Delta(w))\leq \mathrm{p.d.}(\Delta(w))$.
This implies $\mathrm{p.d.}(\Delta(ws))\leq \mathrm{p.d.}(\Delta(w))+1$
and the second statement of Proposition~\ref{p3} follows by induction
from the extreme cases $w=e$ and $w=w_0$ for which it is already 
established.
}
\end{remark}

Now we move to the case of simple modules. 

\begin{proposition}\label{p6}
$\Ext_{\mathcal{O}}^{2l(w_0)-l(w)}(L(w),L(e))\neq 0$, in particular,
$\mathrm{p.d.}(L(w))=2l(w_0)-l(w)$.
\end{proposition}

\begin{proof}
Again the second statement follows from the first statement 
and Proposition~\ref{p1}\eqref{p1.2}. Since $L(w_0)=\Delta(w_0)$,
in the case $w=w_0$ the first statement follows
from Proposition~\ref{p3}. Now we use the inverse induction on $l(w)$. 
Let $s\in W$ be a simple reflection such that $l(sw)<l(w)$. 
Let $m=2l(w_0)-l(w)$. Using the results 
of \cite{AS} we have:
\begin{displaymath}
\begin{array}{lcl}
\mathrm{Ext}_{\mathcal{O}}^{m+1}(\mathrm{T}_sL(w),L(e)) & = & \\
\mathrm{Hom}_{\mathcal{D}^b(\mathcal{O})}(\mathrm{T}_sL(w),L(e)[m+1]) & = & 
\text{(by \cite[Theorems~2.2 and 6.1]{AS})}\\
\mathrm{Hom}_{\mathcal{D}^b(\mathcal{O})}
(\mathcal{L}\mathrm{T}_sL(w),L(e)[m+1]) & = & 
\text{(by \cite[Corollary~4.2]{AS})}\\
\mathrm{Hom}_{\mathcal{D}^b(\mathcal{O})}
(L(w),\mathcal{R}\mathrm{G}_sL(e)[m+1]) & = & 
\text{(by \cite[Corollary~4.2 and 6.2]{AS})}\\
\mathrm{Hom}_{\mathcal{D}^b(\mathcal{O})}
(L(w),L(e)[m]) & = &\\
\mathrm{Ext}_{\mathcal{O}}^{m}(L(w),L(e)).
\end{array}
\end{displaymath}
From the inductive assumption we thus get
$\mathrm{Ext}_{\mathcal{O}}^{m+1}(\mathrm{T}_sL(w),L(e))\neq 0$.
From \cite[Lemma~2.1(3)]{AS} and the right exactness of $\mathrm{T}_s$
it follows that all composition subquotients of 
$\mathrm{T}_sL(w)$ are either of the
form $L(sw)$ or of the form $L(w')$, where $l(w')>l(sw)$. 
From the inductive assumption we have 
$\mathrm{p.d.}(L(w'))\leq m<\mathrm{p.d.}(X)$, which implies
$\mathrm{p.d.}(L(sw))=\mathrm{p.d.}(X)=m+1$. This completes the proof.
\end{proof}

\begin{remark}\label{r7}
{\em Another way to prove the second statement of 
Proposition~\ref{p6} is to use \cite[Theorem~18]{So}, reducing the
question to the Loewy length of some projective module in 
$\mathcal{O}_0$. This Loewy length can then be estimated
using the results from \cite{Ir}. 
}
\end{remark}

\subsection{Costandard modules}\label{s2.2}

An easy corollary from Proposition~\ref{p6} is the following
formula for projective dimensions of costandard modules:

\begin{proposition}\label{p7}
$\mathrm{p.d.}(\nabla(w))=2l(w_0)-l(w)$.
\end{proposition}

\begin{proof}
For $w=w_0$ we have $\nabla(w_0)=L(w_0)$ and the statement follows
from Proposition~\ref{p6}. Now we use the inverse induction on
$l(w)$. Let $s$ be a simple reflection such that $l(ws)<l(w)$. Then
we have the short exact sequence
$\nabla(w)\hookrightarrow\theta_s\nabla(w)\tto \nabla(ws)$.
Since $\theta_s$ is exact and preserves projectives, we have
$\mathrm{p.d.}(\theta_s\nabla(w))\leq \mathrm{p.d.}(\nabla(w))$,
which implies $\mathrm{p.d.}(\nabla(ws))\leq \mathrm{p.d.}(\nabla(w))+1=
2l(w_0)-l(ws)$. On the other hand, for the short exact sequence
$L(ws)\hookrightarrow\nabla(ws)\tto X$ we have that all
simple subquotients of $X$ have the form $L(w')$, where $l(w')>l(ws)$.
Hence, by the inductive assumption, we have $\mathrm{p.d.}(X)<2l(w_0)-l(ws)$,
which implies that $\mathrm{p.d.}(\nabla(ws))=\mathrm{p.d.}(L(ws))$.
The claim follows.
\end{proof}

\begin{remark}\label{r8}
{\em Another way to prove Proposition~\ref{p7} is to use
twisting functors and the results of \cite{AS}, analogously to 
the proofs of Proposition~\ref{p3} and Proposition~\ref{p6}.
}
\end{remark}

\begin{remark}\label{r9}
{\em It is worth mentioning that all the results so far
are obtained without using the Kazhdan-Lusztig conjecture 
(=Theorem).
}
\end{remark}

\subsection{Injective and tilting modules}\label{s2.3}

We are now left to consider the cases of injective and tilting  
modules. It turns out that these are by far more complicated 
than the others. Firstly, we will be forced to use
the Kazhdan-Lusztig conjecture. Secondly, we will not be able to
obtain a description so explicit as above in all cases, and even
in the cases when an explicit description is obtained, the
result is formulated in terms of Kazhdan-Lusztig's combinatorics.
To shorten our notation for $w\in W$ we set
\begin{displaymath}
\mathfrak{t}(w):=\mathrm{p.d.}(T(w)),\quad
\mathfrak{i}(w):=\mathrm{p.d.}(I(w)).
\end{displaymath}
Our main observation about $\mathfrak{t}(w)$ and
$\mathfrak{i}(w)$ is the following:

\begin{theorem}\label{t9}
\begin{enumerate}[(a)]
\item\label{t9.1} Both, $\mathfrak{t}$ and $\mathfrak{i}$, are
constant on the right cells of $W$.
\item\label{t9.2} Both, $\mathfrak{t}$ and $\mathfrak{i}$, are
constant on the left cells of $W$.
\item\label{t9.3} Both, $\mathfrak{t}$ and $\mathfrak{i}$, are
constant on the two-sided cells of $W$.
\end{enumerate}
\end{theorem}

\begin{proof}
The statement \eqref{t9.3} follows immediately from 
\eqref{t9.1} and \eqref{t9.2}.

\begin{proof}[Proof of the statement~\eqref{t9.1}.]
As a consequence of the Kazhdan-Lusztig conjecture, for $w\in W$ 
and a simple reflection, $s\in W$, we have 
(see e.g. \cite[Corollary~5.2.4]{Ir2}):
\begin{equation}\label{eq9-1}
\theta_s\theta_w=
\begin{cases}
\theta_w\oplus\theta_w,& \text{if } ws<w;\\
\theta_{ws}\oplus
\bigoplus_{y<w,ys<y}\mu(y,w)\theta_y,
& \text{if } ws>w,
\end{cases}
\end{equation}
where $\mu(y,w)$ is Kazhdan-Lusztig's $\mu$-function (see
\cite[2.1]{Ir2} or \cite{KL}).

By \cite[Theorem~3.3]{BG} we have $\theta_w P(e)\cong P(w)$
and hence $\theta_w I(e)\cong I(w)$ since $\theta_w$ obviously
commutes with $\star$. Now let $w\in W$ and a simple reflection
$s\in W$ be such that $ws>w$. Since $\theta_s$ is exact and
sends projectives to projectives, applying $\theta_s$ to the
projective resolution of $I(w)=\theta_w I(e)$ and using \eqref{eq9-1}
we obtain that $\mathfrak{i}(ws)\leq \mathfrak{i}(w)$ and
$\mathfrak{i}(y)\leq \mathfrak{i}(w)$ for all
$y$ such that $y<w$, $ys<y$ and $\mu(y,w)\neq 0$. 
In particular, it follows that $\mathfrak{i}$ is monotone
with respect to the right pre-order on $W$ (see e.g. \cite[6.2]{BB}
for details) and thus $\mathfrak{i}$ must be constant on the right cells.

Since $x\mapsto w_0x$ is a bijection on the right cells (see e.g.
\cite[Corollary~6.2.10]{BB}), we have that for  $\mathfrak{t}$ the 
arguments are just the same as for $\mathfrak{i}$, as soon 
as one makes the obvious observation that $\theta_w T(w_0)\cong T(w_0w)$.
\end{proof}

\begin{proof}[Proof of the statement~\eqref{t9.2}.]
The statement~\eqref{t9.2}  is the ``left hand-side version'' of
the statement~\eqref{t9.1}. We would like to prove it using 
analogous arguments, however, for this we will need a
``right hand-side version'' of \eqref{eq9-1}.

\begin{lemma}\label{l10}  
\begin{equation}\label{eq9-2}
\theta_w\theta_s=
\begin{cases}
\theta_w\oplus\theta_w,& \text{if } sw<w;\\
\theta_{sw}\oplus
\bigoplus_{y<w,sy<y}\mu(y,w)\theta_y,
& \text{if } sw>w,
\end{cases}
\end{equation}
\end{lemma}

\begin{proof}
Let $\mathscr{H}$ denote the Hecke algebra of $W$ equipped 
with the standard basis $(H_w)_{w\in W}$. Then there is a unique
antiautomorphism $\sigma$ of $\mathscr{H}$ satisfying
$\sigma(H_s)=H_s$ for any simple reflection $s$. Now
\eqref{eq9-2} is obtained from \eqref{eq9-1} by applying $\sigma$.
\end{proof}

Let $s\in W$ be a simple reflection and $w\in W$. 
Applying $\theta_w$ to the 
short exact sequence $\Delta(sw_0)\hookrightarrow T(sw_0)\tto
\Delta(w_0)$ and observing that $\Delta(sw_0)=\mathrm{G_s}\Delta(w_0)$
(the dual of \cite[(2.3)]{AS}) and $\mathrm{G_s}\theta_w=
\theta_w\mathrm{G_s}$ (the dual of \cite[Theorem~3.2]{AS}), we
get
\begin{equation}\label{eq9-3}
\mathrm{G_s}T(w_0w)\hookrightarrow\theta_w\theta_s T(w)\tto T(w_0w).
\end{equation}

We claim that $\mathrm{p.d.}(\mathrm{G_s}T(w_0w))\leq
\mathrm{p.d.}(T(w_0w))$. Indeed, let us denote $\mathrm{p.d.}(T(w_0w))=m$.
Then for all $i>m$ we have
\begin{displaymath}
\begin{array}{lcl}
\mathrm{Ext}_{\mathcal{O}}^{i}(\mathrm{G_s}T(w_0w),L) & = & \\
\mathrm{Hom}_{\mathcal{D}^b(\mathcal{O})}(\mathrm{G_s}T(w_0w),L[i]) & = & 
\text{(by \cite[Theorems~2.2]{AS})}\\
\mathrm{Hom}_{\mathcal{D}^b(\mathcal{O})}
(\mathcal{R}\mathrm{G_s}T(w_0w),L[i]) & = & 
\text{(by \cite[Corollary~4.2]{AS})}\\
\mathrm{Hom}_{\mathcal{D}^b(\mathcal{O})}
(T(w_0w),\mathcal{L}\mathrm{T_s}L[i]).  
\end{array}
\end{displaymath}
The length of a minimal projective resolution of 
$T(w_0w)$ is $m$. By \cite[Theorems~2.1]{AS}, the non-zero homology of
$\mathcal{L}\mathrm{T_s}L[i]$ can occur only in positions $-i$ or $-i-1$.
Since $i>m$ it follows from \cite[Chapter~III(2), Lemma~2.1]{Ha} that 
$\mathrm{Hom}_{\mathcal{D}^b(\mathcal{O})}
(T(w_0w),\mathcal{L}\mathrm{T_s}L[i])=0$ and thus
$\mathrm{p.d.}(\mathrm{G_s}T(w_0w))\leq
\mathrm{p.d.}(T(w_0w))$.

From the previous paragraph and the short exact sequence
\eqref{eq9-3} we derive the inequality
$\mathrm{p.d.}(\theta_w\theta_s T(w_0))\leq\mathrm{p.d.}(T(w_0w))$.
Now from \eqref{eq9-2} it follows that 
$\mathrm{p.d.}(T(w_0y))\leq \mathrm{p.d.}(T(w_0w))$ for 
each $y$ such that $y<w$, $sy<y$ 
such that $\mu(y,w)\neq 0$. In particular, it follows 
that $\mathfrak{t}$ is monotone with respect to the left pre-order on $W$ 
(see e.g. \cite[6.2]{BB} for details) and thus $\mathfrak{t}$ must be 
constant on the left cells. Again, for
$\mathfrak{i}$ the proof is analogous.
\end{proof}
\end{proof}

\begin{example}\label{ex10}
{\rm
If $\mathfrak{g}$ is of type $A_2$, we have
$W=\{e,s,t,st,ts,sts=tst\}$ with the following decomposition into
two-sided cells: $\{e\}\cup\{s,t,st,ts\}\cup\{sts\}$. 
One easily computes the following table
of values for $\mathfrak{t}$ and $\mathfrak{i}$:
\begin{displaymath}
\begin{array}{|r||c|c|c|c|c|c|}
\hline
w & e & s & t & st & ts & sts\\
\hline
\mathfrak{t}(w) & 0 & 1 & 1 & 1 & 1 & 3\\
\hline
\mathfrak{i}(w) & 6 & 2 & 2 & 2 & 2 & 0\\
\hline
\end{array}
\end{displaymath}
}
\end{example}

There is a well-know integral function on $W$,
constant on two-sided cells, namely  Lusztig's function 
$\mathbf{a}:W\to\mathbb{Z}$, defined in \cite{Lu}. If $w\in W$ is an
involution, then $\mathbf{a}(w)=l(w)-2\delta(w)$, where
$\delta(w)$ is the degree of the Kazhdan-Lusztig polynomial 
$P_{1,w}$, which, together with the property of being constant
on two-sided cells,  completely determines $\mathbf{a}$, 
since every two-sides cell contains a (distinguished) involution, 
see \cite{Lu,Lu2} for details. In particular, if $W_S$ is a 
parabolic subgroup of $W$ and $w_0^S$ is the longest element in 
$W_S$, we have $\mathbf{a}(w_0^S)=l(w_0^S)$. Comparing the values of
$\mathbf{a}$ with Example~\ref{ex10} and other examples 
leads to the following conjecture:

\begin{conjecture}\label{con11}
For all $w\in W$ we have
\begin{enumerate}[(a)]
\item\label{con11.1} $\mathfrak{t}(w)=\mathbf{a}(w)$; 
\item\label{con11.2} $\mathfrak{i}(w)=2\mathbf{a}(w_0w)$. 
\end{enumerate}
\end{conjecture}

\begin{theorem}\label{t12}
Conjecture~\ref{con11} is true if $\mathfrak{g}=\mathfrak{sl}_n$.
\end{theorem}

\begin{proof}
We start by proving Conjecture~\ref{con11}\eqref{con11.1}.

First we observe that in the case $\mathfrak{g}=\mathfrak{sl}_n$ 
every two-sided cell of $W$
contains an element of the form $w_0^S$, where $W_S$ is a parabolic
subgroup of $W$. Indeed, from \cite[Theorem~6.5.1]{BB} we have 
that there is a  bijection between the two-sided cells of 
$\mathcal{S}_n$ and partitions
of $n$. Using \cite[Theorem~6.5.1]{BB} and
\cite[Theorem~3.6.6]{Sa} one gets that the two-sided cell of
$W\cong \mathcal{S}_n$, corresponding to the partition $\lambda\vdash n$,
consists of all $w\in \mathcal{S}_n$, which correspond to standard tableaux
of shape $\lambda$ via the Robinson-Schensted correspondence.
Now if $w_0^S$ is the longest element in some parabolic subgroup of
type $\lambda$, a direct calculation shows that 
the Robinson-Schensted correspondence associates with $w_0^S$ the
partition, which is conjugate to $\lambda$. As a corollary we get that 
every  two-sided cell indeed contains some $w_0^S$. 

Fix now some two-sided cell, say $\mathcal{C}$, and assume that 
it contains $w_0^S$ for some $S$. Because of the properties of
$\mathbf{a}$, listed above,  Conjecture~\ref{con11}\eqref{con11.1}
would follow if we would prove that $\mathrm{p.d.}(T(w_0^S))=l(w_0^S)$.
Assume further that $W_S$ corresponds  to the partition $\lambda$. 
From  \cite[Theorem~6.2.10]{BB}  and \cite[Theorem~3.2.3]{Sa} we get 
that $\mathcal{C}$ also contains an element of the form $w_0w_0^{S'}$, 
where $S'$ corresponds to the conjugate $\lambda'$ of $\lambda$.

Decompose $\theta_{w_0^{S'}}=\theta_{w_0^{S'}}^{out}\theta_{w_0^{S'}}^{on}$,
where $\theta_{w_0^{S'}}^{on}$ is the translation onto the ``most singular''
$S'$-wall, and $\theta_{w_0^{S'}}^{out}$ is the translation out of this
wall. Let further the $w_0^{S'}$-singular block $\mathcal{O}_{\mu}$ be the
image of $\theta_{w_0^{S'}}^{on}$, applied to $\mathcal{O}_0$. Finally,
let $X$ denote the simple Verma module in $\mathcal{O}_{\mu}$.
Then $\theta_{w_0^{S'}}^{on} T(w_0w_0^{S'})\cong X^{\oplus|W_{S'}|}$
and $\theta_{w_0^{S'}}^{out}X\cong T(w_0w_0^{S'})$. Since
translation functors are exact and preserve projectives, we get
$\mathrm{p.d.}(T(w_0w_0^{S'}))=\mathrm{p.d.}(X)$.

The Koszul dual of $\mathcal{O}_{\mu}$ is the regular block of the
$S'$-parabolic category $\mathcal{O}^{\mathfrak{p}}$, see
\cite[Theorem~3.10.2]{BGS}. In particular, via the Koszul duality
$\mathrm{p.d.}(X)$ becomes equal to $m-1$, where $m$ is
the Loewy length of the projective generalized Verma module 
in $\mathcal{O}^{\mathfrak{p}}_0$.
By \cite[Corollary~3.1]{IS}, since $w_0^{S'}$ corresponds to the
partition $\lambda'$, $m-1$ is equal to 
length of the longest element in some parabolic subgroup of
$W$ corresponding to the partition conjugate to $\lambda'$, that
is  to $\lambda$. We finally get that 
\begin{displaymath}
\mathfrak{t}(w_0^S)=\mathfrak{t}(w_0w_0^{S'})=
\mathrm{p.d.}(X)=l(w_0^S).
\end{displaymath}

Now we prove  Conjecture~\ref{con11}\eqref{con11.2} using Conjecture~\ref{con11}\eqref{con11.1}.
In fact, after Conjecture~\ref{con11}\eqref{con11.1} is proved,
one has only to show that  
$\mathfrak{i}(w_0^S)=2\mathfrak{t}(w_0w_0^S)$.
We again decompose 
$\theta_{w_0^{S}}=\theta_{w_0^{S}}^{out}\theta_{w_0^{S}}^{on}$.
We have the singular simple Verma module $X$ such that
$\theta_{w_0^{S}}^{out}X\cong T(w_0w_0^S)$
(and $\theta_{w_0^{S}}^{on}T(w_0w_0^S)\cong X^{\oplus |W_S|}$). 
We also have the singular dominant dual Verma module $Y$ 
such that $\theta_{w_0^{S}}^{out}Y\cong I(w_0^S)$
(and $\theta_{w_0^{S}}^{on}I(w_0^S)\cong Y^{\oplus |W_S|}$). In 
particular, we have $\mathrm{p.d.}(T(w_0w_0^S))=
\mathrm{p.d.}(X)=m$ and $\mathrm{p.d.}(I(w_0^S))=
\mathrm{p.d.}(Y)=n$. So we have to show that $n=2m$.
Taking the Koszul dual we get that $m+1$ equals the Loewy
length of the projective standard module in some regular
block of the parabolic category $\mathcal{O}^{\mathfrak{p}}$.

Let $Z$ denote the simple socle of  $Y$. Then the projective 
dimension of $Z$ equals, via Koszul duality, to $x-1$, where $x$ 
is the  Loewy length of some projective-injective module in
$\mathcal{O}^{\mathfrak{p}}_0$. By \cite[Theorem~5.2(1)]{MS}, 
all projective-injective modules in
$\mathcal{O}^{\mathfrak{p}}_0$ have the same Loewy length.
By \cite[Theorem~5.2(2)]{MS}, the projective generator of
$\mathcal{O}^{\mathfrak{p}}_0$ is a submodule of a projective-injective
module in $\mathcal{O}^{\mathfrak{p}}_0$. It follows that 
projective-injective modules in $\mathcal{O}^{\mathfrak{p}}_0$ have
the maximal possible Loewy length. Thus $\mathrm{p.d.}(Z)$
equals the global dimension of $\mathcal{O}^{\mathfrak{p}}_0$. Since
$Z$ is in the socle of $Y$ and has the maximal possible
projective dimension, from the long exact sequence
in homology it follows that $n=\mathrm{p.d.}(Y)=
\mathrm{p.d.}(Z)=x-1$. Now $n=x-1=2m$ follows from
\cite[Corollary~3.1]{IS}. This completes the proof.
\end{proof}

\begin{remark}\label{rem14}
{\rm
The main difficulty to extend the above arguments to the case 
of arbitrary $\mathfrak{g}$ seems to be the fact that, in general, 
not every two-sided cell contains some element of the form $w_0^S$. 
In fact, Jian-yi Shi has informed me that  in type $D_4$ 
some two-sided cell with $\mathbf{a}$-value $7$ does not
contain any such element. I have no idea how to estimate the values of
$\mathfrak{t}$ and $\mathfrak{i}$ on elements of such cells. In the
general case I can not even prove that $\mathfrak{t}(s)=1$ for a 
simple reflection $s\in W$.
}
\end{remark}

\begin{remark}\label{rem15}
{\rm
The functor $\mathrm{T}=\mathrm{T}_{w_0}$ is exactly the version of 
Arkhipov's functor used in \cite{So2} to establish Ringel's
self-duality of $\mathcal{O}$. In particular, 
$\mathrm{T}P(w)\cong T(w_0w)$ for all $w\in W$.
Using \cite[Corollary~4.2]{AS}, for every $w\in W$ and 
$i\in \mathbb{Z}$ we have
\begin{multline*}
\mathrm{Hom}_{\mathcal{D}^b(\mathcal{O})}(T(w_0w),L[i])=
\mathrm{Hom}_{\mathcal{D}^b(\mathcal{O})}(\mathcal{L}\mathrm{T}P(w),L[i])=\\=
\mathrm{Hom}_{\mathcal{D}^b(\mathcal{O})}(P(w),\mathcal{R}\mathrm{G}L[i]).
\end{multline*}
This shows that Conjecture~\ref{con11} is closely connected
to the understanding of $\mathcal{R}\mathrm{G}$ applied
to simple modules, that
is to the understanding of the homology of the complex
$\mathrm{G}\mathcal{I}^{\bullet}$, where 
$\mathcal{I}^{\bullet}$ is an injective resolution of $L$. 
We remark that $\mathcal{I}^{\bullet}$ is a projective object in 
$\mathscr{LC}(\mathtt{I})$; and  $\mathrm{G}\mathcal{I}^{\bullet}$ 
is a projective object in the  category $\mathscr{LC}(\mathtt{T})$
(see \cite[Proposition~11]{MOS}).
These categories will appear later on in the paper, where we will
also try study the connection mentioned above in more details.
}
\end{remark}

\subsection{Shuffled Verma modules}\label{s2.4}

There is a very special class of modules in $\mathcal{O}_0$,
called {\em shuffled Verma modules}, which were introduced in
\cite{Ir3} as modules, corresponding to the principal series
modules. Using \cite[Section~3]{AL}
for $x,y\in W$ we define the corresponding shuffled Verma
module 
\begin{displaymath}
\Delta(x,y)=\mathrm{T}_x\Delta(y)
\end{displaymath}
(as these modules are defined using the twisting functors,
sometimes they are also called {\em twisted Verma modules},
however, we will use the name {\em shuffled Verma modules}
as in the original paper \cite{Ir3}).
In particular, using \cite[(2.3) and Theorem~2.3]{AS} for any 
$w\in W$ we have 
\begin{displaymath}
\begin{array}{rclcccrcl}
\Delta(e,w)& \cong & \Delta(w),&&&&
\Delta(w,w_0)& \cong & \nabla(ww_0),\\
\Delta(w,e)& \cong & \Delta(w),&&&&
\Delta(w_0,w)& \cong & \nabla(w_0w).
\end{array}
\end{displaymath}

For shuffled Verma modules we have the following statement, which
includes Proposition~\ref{p3} and Proposition~\ref{p7} as special cases:

\begin{proposition}\label{p15}
For $x,y\in W$ we have $\mathrm{p.d.}(\Delta(x,y))=l(x)+l(y)$.
\end{proposition}

\begin{proof}
First let us prove that $\mathrm{p.d.}(\Delta(x,y))\leq l(x)+l(y)$
by induction on $l(x)$. If $x=e$, the statement follows from 
Proposition~\ref{p7}. Let now $x=sz$, where $s$ is a simple reflection
and $l(z)<l(x)$. Since $\Delta(x,y)=\mathrm{T}_s\Delta(z,y)$,
for $i>l(x)+l(y)$ we have
\begin{equation}\label{eqp15-1}
\begin{array}{lcl}
\mathrm{Ext}_{\mathcal{O}}^{i}(\mathrm{T}_s\Delta(z,y),L) & = & \\
\mathrm{Hom}_{\mathcal{D}^b(\mathcal{O})}(\mathrm{T}_s\Delta(z,y),L[i]) & = & 
\text{(by \cite[Theorems~2.2]{AS})}\\
\mathrm{Hom}_{\mathcal{D}^b(\mathcal{O})}
(\mathcal{L}\mathrm{T}_s\Delta(z,y),L[i]) & = & 
\text{(by \cite[Corollary~4.2]{AS})}\\
\mathrm{Hom}_{\mathcal{D}^b(\mathcal{O})}
(\Delta(z,y),\mathcal{R}\mathrm{G_s}L[i]).  
\end{array}
\end{equation}
By the assumption of induction we know that the projective resolution
of $\Delta(z,y)$ has length at most $l(x)+l(y)-1$. By the dual of 
\cite[Theorems~2.2]{AS}, non-zero homology of $\mathcal{R}\mathrm{G_s}L[i]$
can occur only in positions $-i,-i+1<-(l(x)+l(y)-1)$. Hence,
using \cite[Chapter~III(2), Lemma~2.1]{Ha}, we get that
$\mathrm{Hom}_{\mathcal{D}^b(\mathcal{O})}
(\Delta(z,y),\mathcal{R}\mathrm{G_s}L[i])=0$.

Now it is enough to observe that
$\mathrm{Ext}_{\mathcal{O}}^{l(x)+l(y)}(\Delta(x,y),L(e))\neq 0$.
We use induction on $l(x)+l(y)$. If $l(x)=0$, this is proved in
Proposition~\ref{p3}.  If $l(x)>1$ this follows from 
the inductive assumption and \eqref{eqp15-1}
using \cite[Corollary~2.2]{AS}. This completes the proof.
\end{proof}

\begin{remark}\label{rem21-7}
{\rm
Twisted tilting modules $\mathrm{T}_x T(y)$, $x,y\in W$, were
studied in \cite{St2}. One can also consider the twisted projective
modules $\mathrm{T}_x P(y)$, $x,y\in W$ (for $x=w_0$ the latter
coincide with the usual tilting modules). It is a natural question 
to determine the projective dimension of these modules. However, this
question seems to be even more complicated than the corresponding
question for the usual tilting modules. The main reason is that,
in contrast to the usual tilting modules, for twisted
tilting or twisted projective modules the function of projective
dimension will be constant only on the appropriate right cells, but 
not on the two-sided cells in the general case.
}
\end{remark}

\section{On the extension algebra of standard modules}\label{s3}

\subsection{Setup for Koszul quasi-hereditary algebras}\label{s3.1}

Let $\Bbbk$ be an algebraically closed field.
Let $\mathtt{A}=\oplus_{i\in\mathbb{Z}}\mathtt{A}_i$ be a 
{\em positively graded} $\Bbbk$-algebra, that is 
$\dim \mathtt{A}_i=0$ for all $i<0$; $\dim \mathtt{A}_i<\infty$
for all $i$; and $\mathtt{A}_0=\oplus_{\lambda\in\Lambda}\Bbbk
e_{\lambda}$, where $1=\sum_{\lambda\in\Lambda}e_{\lambda}$ is a
fixed decomposition of $1$ into a sum of pairwise orthogonal 
primitive idempotents. We denote by $\mathtt{A}^!$ the
{\em quadratic dual} of $\mathtt{A}$, see e.g. 
\cite[Section~6]{MO2}.

Let $\mathtt{A}\mathrm{-fgmod}$ denote the category of all
graded $\mathtt{A}$-modules with finite-dimensional graded components.
Morphisms in this category are homogeneous maps of degree $0$ between
graded modules. Under our assumptions, this category contains several
natural classes of modules. To each $\lambda\in \Lambda$ 
there correspond the graded projective module  
$\mathtt{P}(\lambda)=\mathtt{A}e_{\lambda}$, 
its simple quotient $\mathtt{S}(\lambda)$,
and the injective hull $\mathtt{I}(\lambda)$ of $\mathtt{S}(\lambda)$.
Assume further that $\mathtt{A}$ is quasi-hereditary with respect
to some order $\leq$ on $\Lambda$. Then we also have the corresponding
graded standard module $\mathtt{\Delta}(\lambda)$, the
graded costandard module $\mathtt{\nabla}(\lambda)$,  and the graded 
tilting modules $\mathtt{T}(\lambda)$, (see for example \cite{Zh}).
As before we set $\mathtt{P}=\oplus_{\lambda\in\Lambda}\mathtt{P}(\lambda)$
and analogously for all other types of modules.
We have that the canonical surjections 
$\mathtt{P}(\lambda)\tto \mathtt{\Delta}(\lambda)\tto \mathtt{S}(\lambda)$
and $\mathtt{T}(\lambda)\tto \mathtt{\nabla}(\lambda)$,  and the
canonical injections $\mathtt{S}(\lambda)\hookrightarrow 
\mathtt{\nabla}(\lambda)\hookrightarrow \mathtt{I}(\lambda)$ and
$\mathtt{\Delta}(\lambda)\hookrightarrow \mathtt{T}(\lambda)$
are morphisms in $\mathtt{A}\mathrm{-fgmod}$. As before
$\langle k\rangle$ denotes the shift of grading.

Denote by $\mathscr{LC}(\mathtt{P})$ (resp. 
$\mathscr{LC}(\mathtt{T})$) the category, 
whose objects are all complexes $\mathcal{X}^{\bullet}$ such that 
$\mathcal{X}^{i}\in \mathrm{add}(\mathtt{P}\langle i\rangle)$
(resp. $\mathrm{add}(\mathtt{T}\langle i\rangle)$) for all $i$,
and morphisms are all morphisms of complexes. From the positivity
of the grading it follows that the only homotopy between two objects
of $\mathscr{LC}(\mathtt{P})$ is the trivial one. The grading on $\mathtt{A}$
automatically induces a grading on the {\em Ringel dual}
$R(\mathtt{A})=\mathrm{End}_{\mathtt{A}}(\mathtt{P})^{\mathrm{op}}$. If this
grading is positive (which is not true in general), then the only 
homotopy between two objects of $\mathscr{LC}(\mathtt{T})$ is the trivial one
(see \cite[Section~6]{MO2}).  The category $\mathscr{LC}(\mathtt{P})$ is 
equivalent to $\mathtt{A}^!\mathrm{-fgmod}$ and the category
$\mathscr{LC}(\mathtt{T})$ is  equivalent to $R(\mathtt{A})^!\mathrm{-fgmod}$,
see e.g. \cite[Section~6]{MO2}.

Assume now that both the minimal titling coresolution of 
$\mathtt{\Delta}$ and the minimal tilting resolution of 
$\mathtt{\nabla}$ are objects in  $\mathscr{LC}(\mathtt{T})$. In particular,
this implies (see \cite[Theorem~7]{MO2}) that $\mathtt{A}$ is standard
Koszul in the sense of \cite{ADL}. Hence the algebra
$R(\mathtt{A})^!$ is quasi-hereditary. Certainly $R(\mathtt{A})^!$
inherits a grading. Finally, we assume that the induced grading on
$R(R(\mathtt{A})^!)$ is positive (which means that $\mathtt{A}$
is {\em balanced} in the sense of \cite[Section~6]{MO2}).

\subsection{Bigraded extension algebra of standard modules}\label{s3.2}

Consider the full subcategory of 
$\mathcal{D}^b(\mathtt{A}\mathrm{-fgmod})$,
whose objects are $\mathtt{\Delta}(\lambda)\langle i\rangle[j]$, 
where $\lambda\in\Lambda$, $i,j\in\mathbb{Z}$. The group 
$\mathbb{Z}^2$ acts freely on this category by shifting the 
grading and the position in the complex. This induces a canonical $\mathbb{Z}^2$-grading on the (originally {\em ungraded}) 
Yoneda Ext-algebra 
$\mathrm{Ext}^*_{A}(\Delta)$, see e.g. \cite{DM}.
This $\mathbb{Z}^2$-graded
algebra has two natural $\mathbb{Z}$-graded subalgebras. The
first one the $\mathbb{Z}$-graded algebra
$\mathrm{End}^*_{A}(\Delta)$ of all homomorphisms
between graded standard modules obtained  in the folowing way:
Consider the full subcategory of 
$\mathcal{D}^b(\mathtt{A}\mathrm{-fgmod})$,
whose objects are $\mathtt{\Delta}(\lambda)\langle i\rangle$, 
where $\lambda\in\Lambda$, $i\in\mathbb{Z}$. The group
$\mathbb{Z}$ acts freely on this category by shifting the
garding. $\mathrm{End}^*_{A}(\Delta)$ is the $\mathbb{Z}$-graded 
algebra obtained as the quotient of this action. The
second subalgebra is the $\mathbb{Z}$-graded algebra
$\mathrm{Lext}^*_{A}(\Delta)$ of all {\em linear extensions}
defined  in the folowing way: Consider the full subcategory of 
$\mathcal{D}^b(\mathtt{A}\mathrm{-fgmod})$,
whose objects are $\mathtt{\Delta}(\lambda)\langle i\rangle[-i]$, 
where $\lambda\in\Lambda$, $i\in\mathbb{Z}$. The group
$\mathbb{Z}$ acts freely on this category via $\langle i\rangle[-i]$,
$i\in\mathbb{Z}$. $\mathrm{Lext}^*_{A}(\Delta)$ is the 
$\mathbb{Z}$-graded algebra obtained as the quotient of this action.
Our main general result in
this section is the following fairly obvious observation, which,
however, will have some interesting applications to the category
$\mathcal{O}$.

\begin{proposition}\label{p2.01}
Let $\mathtt{A}$ be balanced. Then the Yoneda extension algebras of
standard modules for $\mathtt{A}$ and $R(\mathtt{A})^!$ are
canonically isomorphic as $\mathbb{Z}^2$-graded algebras.
This isomorphism induces the following isomorphisms of
$\mathbb{Z}$-graded subalgebras:
\begin{displaymath}
\begin{array}{rcl}
\mathrm{End}^*_{A}(\Delta) & \cong &
\mathrm{Lext}^*_{R(\mathtt{A})^!}(\Delta),\\
\mathrm{Lext}^*_{A}(\Delta) & \cong &
\mathrm{End}^*_{R(\mathtt{A})^!}(\Delta).
\end{array}
\end{displaymath}
\end{proposition}

\begin{proof}
Since $\mathtt{A}$ is balanced, then both $\mathtt{A}$ and 
$R(\mathtt{A})$ are quasi-hereditary and Koszul. The Ringel
and Koszul dualities induce equivalences between the corresponding
bounded derived categories of graded modules. By \cite[Theorem~9]{MO2},
standard modules for  $\mathtt{A}$ and $R(\mathtt{A})^!$
can be identified via these dualities. The first part of the
claim follows. The second part follows from the identification of
standard modules, given in \cite[Theorem~9]{MO2}.
\end{proof}

\subsection{Applications to the category $\mathcal{O}$}\label{s3.3}

Proposition~\ref{p2.01} can immediately be applied to the graded
algebra $\mathtt{A}$ of the principal block of the category 
$\mathcal{O}$. Namely, in the 
notation of Section~\ref{s1} we have.

\begin{theorem}\label{t2.02}
\begin{enumerate}[(a)]
\item\label{t2.02.1} There is a {\em non-trivial} automorphism of the
$\mathbb{Z}^2$-graded algebra $\mathrm{End}^*_{A}(\Delta)$,
which swaps $\mathrm{End}^*_{A}(\Delta)$ and  
$\mathrm{Lext}^*_{A}(\Delta)$. In particular, the 
$\mathbb{Z}$-graded algebras $\mathrm{End}^*_{A}(\Delta)$ and  
$\mathrm{Lext}^*_{A}(\Delta)$ are isomorphic.
\item\label{t2.02.2}
$\mathrm{Ext}^i_{\mathtt{A}}(\Delta(x),
\Delta(y)\langle j\rangle)\cong
\mathrm{Ext}^{i+j}_{\mathtt{A}}(\Delta(w_0y^{-1}w_0),
\Delta(w_0x^{-1}w_0)\langle -j\rangle)$ for all elements
$x,y\in W$.
\end{enumerate}
\end{theorem}

\begin{proof}
$\mathtt{A}$ is both Koszul self-dual (\cite[Theorem~18]{So})
and Ringel self dual (\cite[Corollary~2.3]{So2}). Hence the 
first statement follows directly from Proposition~\ref{p2.01}.
The second statement follows by tracking
the correspondence induced by these self-dualities on 
primitive idempotents and \cite[Theorem~21(ii)]{MOS}. 
\end{proof}

The latter statement has some interesting corollaries. The
first one describes the linear extensions between standard
modules:

\begin{corollary}\label{c2.03}
For $x,y\in W$, we have:
\begin{displaymath}
\mathrm{Ext}^i_{\mathtt{A}}(\Delta(x),
\Delta(y)\langle -i\rangle)\cong
\begin{cases}
\mathbb{C}, & x\geq y\text{ and } l(x)-l(y)=i;\\
0, & \text{otherwise}.
\end{cases}
\end{displaymath}
\end{corollary}

\begin{proof}
Theorem~\ref{t2.02} reduces the statement to the
analogous statement for homomorphisms between 
Verma modules. We know that the positive
grading on $\mathtt{A}$ induces a positive
grading on Verma modules. Furthermore, we also
know when homomorphisms between
Verma modules do exist, and that the homomorphism space
between Verma modules is at most one-dimensional
(see \cite[Section~7]{Di}).
Moreover, all Verma modules have the same simple socle. 
So, to get the explicit formula above one 
has to compare the lengths of their graded
filtrations, which can be done using, for example,
\cite[Section~5]{St}.
\end{proof}

\begin{remark}
{\rm
From Corollary~\ref{c2.03} it follows that the
assertion of \cite[Theorem~6]{MO2} requires some additional
assumptions, for example it is sufficient to make
\cite[Assumptions (I) and (II)]{DM}.
}
\end{remark}

Another corollary is the following result of Carlin
(see \cite[(3.8)]{Ca}):

\begin{corollary}\label{c2.04}
For $x,y\in W$, $x\geq y$, we have
$\mathrm{Ext}^{l(x)-l(y)}_{A}(\Delta(x),\Delta(y))\cong \mathbb{C}$.
\end{corollary}

\begin{proof}
Since $A$ is quasi-hereditary with respect to the Bruhat order on
$W$, the projective modules, occurring at the position 
$l(y)-l(x)$ in the minimal 
(linear) projective resolution of $\Delta(x)$, have indexes
$w$ such that $l(w)\leq l(y)$. At the same time all simple modules
in the radical of $\Delta(y)$ have indexes $u$ such that
$l(u)> l(y)$. Hence any non-zero element in the space
$\mathrm{Ext}^{l(x)-l(y)}_{A}(\Delta(x),\Delta(y))$ must
belong to $\mathrm{Ext}^{l(x)-l(y)}_{\mathtt{A}}
(\mathtt{\Delta}(x),\mathtt{\Delta}(y)\langle l(y)-l(x)\rangle)$. 
Now the statement follows from Corollary~\ref{c2.03}.
\end{proof}

\begin{remark}\label{remnew2.01}
{\em
Using the parabolic-singular Koszul duality from \cite{BGS,Ba2}  
and \cite[Appendix]{MO2} one obtains that the extension
algebras of standard modules for parabolic and corresponding
singular blocks (respectively, pairs of corresponding
pa\-ra\-bo\-lic-\-sin\-gu\-lar blocks) 
are also isomorphic as bigraded algebras.  This isomorphism 
again swaps the subalgebra of homomorphisms
with the subalgebra of linear extensions.
}
\end{remark}

\subsection{Several graded subalgebras of the extension algebra of
standard modules}\label{s3.4}

We continue to study the $\mathbb{Z}^2$-graded extension algebra
$\mathrm{Ext}^*_{A}(\Delta)$ of the block $\mathcal{O}_0$. From
the quasi-heredity of $\mathtt{A}$ we immediately obtain the following 
vanishing condition:
$\mathrm{Ext}^i_{\mathtt{A}}(\Delta,\Delta\langle j\rangle)\neq 0$
implies $i\geq 0$ and $j\geq -i$. It follows that the following
induces a natural {\em positive} $\mathbb{Z}$-grading on 
$\mathtt{E}:=\mathrm{Ext}^*_{A}(\Delta)$ (in the sense of 
\cite[2.1]{MOS}):
\begin{displaymath}
\mathtt{E}_k=\bigoplus_{2i+j=k}
\mathrm{Ext}^i_{A}(\Delta,\Delta\langle j\rangle),
\quad k\in\mathbb{Z}.
\end{displaymath}
In particular, both $\mathrm{End}^*_{A}(\Delta)$ and  
$\mathrm{Lext}^*_{A}(\Delta)$ become $\mathbb{Z}$-graded
subalgebras of $\mathtt{E}$ in the natural way.

\begin{remark}\label{r2.06}
{\rm 
The natural $\mathbb{Z}$-grading on  $\mathtt{E}$ given by the
degree of the extension is {\em not positive} since the zero
component of this grading (the subalgebra of all homomorphisms) 
is not a semi-simple subalgebra in the general case.
}
\end{remark}

Our first result  here is the following Koszulity statement for the 
subalgebra of all homomorphisms.

\begin{proposition}\label{p2.05}
The  algebra $\mathrm{End}^*_{A}(\Delta)$
is Koszul.
\end{proposition}

\begin{proof}
First I claim that, as a $\mathbb{Z}$-graded algebra, the algebra $\mathrm{End}^*_{A}(\Delta)$ is isomorphic to the incidence algebra
of the poset $W$ with respect to $\leq$. Let us describe 
$\mathrm{End}^*_{A}(\Delta)$ via some quiver with relations.
For $x,y\in W$, $x\geq y$, we have a unique up to scalar injection
$\Delta(x)\hookrightarrow\Delta(y)$. In particular, we can
identify each $\Delta(w)$, $w\in W$, with the corresponding
submodule of $\Delta(e)$. For each $w\in W$ let $v_w$ denote 
some generator of $\Delta(w)$, which we fix.  
If $x,y\in W$, $x\geq y$, let
$\varphi_{x,y}:\Delta(x)\to\Delta(y)$ denote the homomorphism,
such that $\varphi_{x,y}(v_x)=v_x$. Then, by 
\cite[Theorem~7.6.23]{Di}, the arrows in the quiver of
$\mathrm{End}^*_{A}(\Delta)$ are
$\varphi_{x,y}$ such that $x=sy$, where $s$ is a reflection
(not necessarily simple). From the definition of 
$\varphi_{x,y}$ we have that these arrows obviously satisfy 
all relevant commutativity relations. Hence 
$\mathrm{End}^*_{A}(\Delta)$ is a quotient of the 
incidence algebra of the poset $(W,\geq)$. It follows that
the two algebras coincide because they obviously have 
the same dimension.

Now we recall that the M{\"o}bius function of the poset 
$(W,\geq)$ was determined in \cite{Ve}. It equals
$(-1)^{l(x)-l(y)}$ for $x\geq y$. Hence, the Koszulity
of the corresponding incidence algebra follows
from \cite[Theorem~1]{Yu}. This completes the proof.
\end{proof}

In \cite{DM} it is shown that in the multiplicity-free cases
the $\mathbb{Z}$-graded algebra $\mathtt{E}$ is Koszul with 
respect to the positive grading introduced above. This and
Proposition~\ref{p2.05} motivate the following conjecture:

\begin{conjecture}\label{con2.06}
The subalgebra of $\mathtt{E}$ generated by 
$\mathtt{E}_0$ and $\mathtt{E}_1$ is Koszul.
\end{conjecture}

\begin{remark}\label{rem2.07}
{\rm
I do not know if $\mathtt{E}_0$ and $\mathtt{E}_1$ generate
the whole $\mathtt{E}$ in general. I do believe that
they do not. Elements from $\mathtt{E}_1$ correspond
to ``naive'' extensions, which do not take into account
the multiplicities given by the Kazhdan-Lusztig combinatorics.
Additionally, the numerical structure of extensions between 
Verma modules seems to be really complicated, see \cite{GJ,Ca,Bo}.
}
\end{remark}

\subsection{Some remarks on extensions between
Verma modules}\label{s3.5}

As already mentioned, the description of the algebra
$\mathtt{E}$, and even of the dimensions
$\dim\mathrm{Ext}_{\mathtt{A}}^i(\Delta(x),
\Delta(y)\langle j\rangle)$ seems to be a very complicated
problem, see \cite{GJ,Ca,Bo}. A very easy observation
reduces this problem to the description of certain
properties of the funtor
$\mathcal{L}\mathrm{T}_{x}$:

\begin{proposition}\label{p2.10}
Let $x,y\in W$ and $i,j\in \mathbb{Z}$. Then
\begin{displaymath}
\begin{array}{lll}
\dim\mathrm{Ext}_{\mathtt{A}}^i(\Delta(x),
\Delta(y)\langle j\rangle)&=&
\left[\mathcal{R}^i\mathrm{G}_{x^{-1}}\Delta(y)\langle j\rangle:L(e)
\right]\\ &=&
[\mathcal{L}_i\mathrm{T}_{x^{-1}}\nabla(y)\langle -j\rangle:L(e)].
\end{array}
\end{displaymath}
\end{proposition}

\begin{proof}
Taking into account that twisting functors are
gradable (see \cite[Appendix]{MO2} or \cite[page~28]{FKS}), 
we compute:
\begin{displaymath}
\begin{array}{lcl}
\mathrm{Ext}_{\mathtt{A}}^i(\Delta(x),
\Delta(y)\langle j\rangle) & = & \\
\mathrm{Hom}_{\mathcal{D}^b(\mathtt{A})}(\Delta(x),
\Delta(y)\langle j\rangle[i]) & = & 
\text{(by \cite[(2.3)]{AS})}\\
\mathrm{Hom}_{\mathcal{D}^b(\mathtt{A})}(\mathrm{T}_x\Delta(e),
\Delta(y)\langle j\rangle[i]) & = & 
\text{(by \cite[Theorem~2.2]{AS})}\\
\mathrm{Hom}_{\mathcal{D}^b(\mathtt{A})}(\mathcal{L}\mathrm{T}_x\Delta(e),
\Delta(y)\langle j\rangle[i]) & = & 
\text{(by \cite[Corollary~4.2]{AS})}\\
\mathrm{Hom}_{\mathcal{D}^b(\mathtt{A})}(\Delta(e),
\mathcal{R}\mathrm{G}_{x^{-1}}\Delta(y)\langle j\rangle[i]) & = & 
\text{($\Delta(e)$~-- projective)}
\\
\left[\mathcal{R}^i\mathrm{G}_{x^{-1}}\Delta(y)\langle j\rangle:L(e)
\right]
& = & 
\text{(by duality)}\\
\left[\mathcal{L}_i\mathrm{T}_{x^{-1}}\nabla(y)\langle -j\rangle:L(e)
\right].
 &&
\end{array}
\end{displaymath}
\end{proof}

\begin{remark}
{\rm 
Since both, the twisting and the shuffling functors, 
are auto-equivalences of $\mathcal{D}^b(\mathcal{O}_0)$
(see \cite[Corollary~4.2]{AS} and \cite[Theorem~5.7]{MS2}), we have 
\begin{displaymath}
\begin{array}{lll}
\mathrm{Ext}^i_{\mathtt{A}}(\Delta(sx),\Delta(sy))=
\mathrm{Ext}^i_{\mathtt{A}}(\Delta(x),\Delta(y))
& \text{ if }& sx>x, sy>y;\\ 
\mathrm{Ext}^i_{\mathtt{A}}(\Delta(xs),\Delta(ys))=
\mathrm{Ext}^i_{\mathtt{A}}(\Delta(x),\Delta(y))
& \text{ if }& xs>x, ys>y.\\ 
\end{array}
\end{displaymath}
Since both, the twisting and the shuffling functors, are gradable,
the above formula admits a natural graded analogue.
In many cases, but not in all, this formula can be applied to reduce 
extensions to the case of extensions into the projective standard 
module. In  particular, the latter extensions deserve special attention.
}
\end{remark}

Here we would like to present one application of the above
technique, which gives (from my point of view) 
a fairly unexpected description of the
$\mathrm{Ext}^1$-space into the projective standard module.
For $x\in W$ with a fixed reduced decomposition 
$x=s_1\cdots x_k$ we denote by $\overline{l}(x)$ the number
of different simple reflections occurring in this
reduced decomposition (for example $\overline{l}(sts)=2$
if $s$ and $t$ do not commute). Since any two reduced decompositions
can be obtained from each other by applying braid relations
only, it follows that $\overline{l}(x)$ does not
depend on the reduced decomposition of $x$.

\begin{theorem}\label{t2.12}
\begin{displaymath}
\dim\mathrm{Ext}^1_{\mathtt{A}}
(\mathtt{\Delta}(x),\mathtt{\Delta}(e)\langle j\rangle)=
\begin{cases}
\overline{l}(x), &  \text{ if } j=l(x)-2;\\
0, & \text{ otherwise.}
\end{cases}
\end{displaymath}
\end{theorem}

\begin{proof}
We start with a special case:

\begin{lemma}\label{l2.14}
The statement of Theorem~\ref{t2.12} is true in the case
$x=w_0$.
\end{lemma}

\begin{proof}
Let $\Delta(e)\hookrightarrow X\tto \Delta(w_0)$ be a non-split
extension. Since $\Delta(w_0)$ is simple and $\Delta(e)$ has
simple socle $L(w_0)$ it follows that $X$ has simple socle
$L(w_0)$. In particular, $X\hookrightarrow P(w_0)$. Since both
$\Delta(e)$ and $\Delta(w_0)$ have central characters it follows
that $X$ is annihilated by the second power of the corresponding
maximal ideal of the center. By \cite[Proposition~2.12]{Ba}, this
means that $X$ is a submodule of the submodule $Y\subset P(w_0)$, 
which is uniquely determined via $\Delta(e)\hookrightarrow Y\tto 
\bigoplus_{s:l(s)=1}\Delta(s)$. Since each $\Delta(s)$ has simple
socle $\Delta(w_0)$ and no other occurrences of $\Delta(w_0)$ in the
composition series, we have that $X$ is even a submodule of the
submodule $Z$ of $Y$ such that  $\Delta(e)\hookrightarrow Z\tto 
\Delta(w_0)^{\oplus k }$, where $k=|\{s:l(s)=1\}|$. 
Since $Z$ has simple socle, it follows that 
$\dim\mathrm{Ext}^{1}_A(\Delta(w_0),\Delta(e))$ equals the
number of simple roots, which obviously equals $\overline{l}(w_0)$.
Now the necessary statement follows by tracking the grading
using \cite{St} and \cite[Appendix]{MO2}.
\end{proof}

Now we go to the general case. Our strategy is: we first 
establish a lower bound and then prove that it is in fact the
real value. As it was already done in Lemma~\ref{l2.14}, it is
easier to prove the ungraded version and then just track the necessary
grading using \cite{St} and \cite[Appendix]{MO2}.

Set $y=x^{-1}$ and observe that $\overline{l}(x)=\overline{l}(y)$.
Now consider the the short exact sequence
\begin{equation}\label{eq2.12.1}
0 \to \Delta(w_0)\to P(w_0)\to Coker\to 0.
\end{equation}
Note that  $P(w_0)$ is injective. Let $\alpha$ denote the
natural transformation from $\mathrm{ID}$ to $\mathrm{G}_{y}$
given by \cite[2.3]{KM}. Observe that $\alpha$ is injective on
all modules from \eqref{eq2.12.1} since they all have Verma flags
(this follows, for example, from the dual of
\cite[Proposition~5.4]{AS}).
Further note that $\alpha$ is an isomorphism on both
$\Delta(e)$ and $P(w_0)$ because of the projectivity of these
two modules (by the dual of
\cite[Corollary~9]{KM}).
Now, applying $\mathrm{G}_{y}$ to \eqref{eq2.12.1} yields to the
following commutative diagram with exact columns and rows:
\begin{displaymath}
\xymatrix{
\Delta(e)\ar@{^{(}->}[r]\ar[d]^{\wr} & 
P(w_0)\ar@{->>}[r]\ar[d]^{\wr} & 
Coker\ar@{^{(}->}[d]^{\alpha_{Coker}}\\
\mathrm{G}_{y}\Delta(e)\ar@{^{(}->}[r] & 
\mathrm{G}_{y}P(w_0)\ar[r]^f & 
\mathrm{G}_{y}Coker\ar@{->>}[r]\ar@{->>}[d]
& \mathcal{R}^1\mathrm{G}_{y}\Delta(e)\\
 & & X 
}
\end{displaymath}
From this diagram  we have that the heads of the image of 
both $f$ and $\alpha_{Coker}$ are isomorphic to $L(w_0)$ and that
the multiplicity of $L(w_0)$ in both $X$ 
and $\mathcal{R}^1\mathrm{G}_{y}\Delta(e)$ is $0$. This
implies that the kernel of both $f$ and $\alpha_{Coker}$ is
the trace of $P(w_0)$ in $\mathrm{G}_{y}Coker$, in particular,
$X=\mathcal{R}^1\mathrm{G}_{y}\Delta(e)$.

Let now $Y=\oplus_{s:l(s)=1}\Delta(s)$. Then we have the following 
short exact sequence: $Y\hookrightarrow Coker\tto Coker'$, where again
all modules have Verma flags. Applying $\mathrm{G}_{y}$ and using
the Snake Lemma gives
the following commutative diagram with exact rows and columns:
\begin{equation}\label{eq2.12.2}
\xymatrix{
Y \ar@{^{(}->}[r]\ar@{^{(}->}[d]^{\alpha_Y}& 
Coker\ar@{^{(}->}[d]^{\alpha_{Coker}}\ar@{->>}[r] & 
Coker'\ar@{^{(}->}[d]^{\alpha_{Coker'}} \\
\mathrm{G}_{y}M\ar@{^{(}->}[r]\ar@{->>}[d] & 
\mathrm{G}_{y}Coker\ar@{->>}[d]\ar@{->}[r] & 
\mathrm{G}_{y}Coker' \\
Z\ar@{^{(}->}[r] &
\mathcal{R}^1\mathrm{G}_{y}\Delta(e) 
}
\end{equation}
Let $S_1$ denote the set of all simple roots which appear in a
reduced expression of $y$, and let $S_2$ denote the set of all
other simple roots. From the dual of \cite[Theorem~2.3]{AS} we get
\begin{displaymath}
Z\cong \Delta(e)^{\oplus |S_1|} \oplus
\bigoplus_{s\in S_2} \Delta(s).
\end{displaymath}
In particular, we obtain that $[Z:L(e)]=|S_1|$ and hence 
$[\mathcal{R}^1\mathrm{G}_{y}\Delta(e):L(e)]\geq |S_1|$ because of 
the third row of \eqref{eq2.12.2}. This is our lower bound.

Now to prove that this lower bound gives the exact value, we
write $\mathrm{G}_{w_0}=\mathrm{G}_{z}\mathrm{G}_{y}$,
where $z=w_0x$ and note
that the natural transformation from $\mathrm{ID}$ to
$\mathrm{G}_{w_0}$ can be obviously written as the composition 
of the natural transformation from $\mathrm{ID}$ to
$\mathrm{G}_{y}$ with the natural transformation
from $\mathrm{ID}$ to $\mathrm{G}_{z}$, the latter being
restricted to the
image of $\mathrm{G}_{y}$. This implies that the diagram
\eqref{eq2.12.2} can be extended to the following
commutative diagram with exact rows and columns:
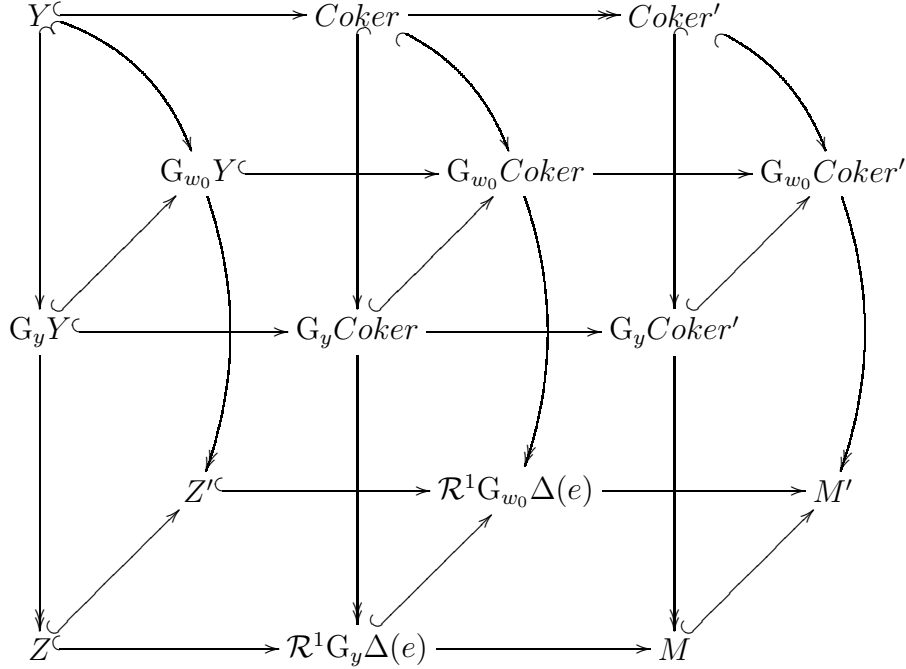
\begin{displaymath}
\xymatrix@!=3pc{
Y \ar@{^{(}->}[rr]\ar@{^{(}->}[dd]\ar@{_{(}->}@/^1pc/[dr]&& 
Coker\ar@{->>}[rr]\ar@{^{(}->}[dd]\ar@{_{(}->}@/^1pc/[dr] && 
Coker'\ar@{^{(}->}[dd]\ar@{_{(}->}@/^1pc/[dr] \\
 & 
 \mathrm{G}_{w_0}Y\ar@{^{(}->}[rr]\ar@{->>}@/^1pc/[dd] && 
 \mathrm{G}_{w_0}Coker\ar@{->>}@/^1pc/[dd]\ar[rr] && 
 \mathrm{G}_{w_0}Coker'\ar@{->>}@/^1pc/[dd] \\
\mathrm{G}_{y}Y \ar@{^{(}->}[rr]\ar@{->>}[dd]\ar@{^{(}->}[ur]&& 
\mathrm{G}_{y}Coker\ar@{->>}[dd]\ar[rr]\ar@{^{(}->}[ur] && 
\mathrm{G}_{y}Coker'\ar@{->>}[dd] \ar@{^{(}->}[ur]\\
 & 
 Z'\ar@{^{(}->}[rr] && 
 \mathcal{R}^1\mathrm{G}_{w_0}\Delta(e)\ar[rr] && 
 M' \\
Z\ar@{^{(}->}[rr]\ar@{^{(}->}[ur]&& 
\mathcal{R}^1\mathrm{G}_{y}\Delta(e)\ar[rr]\ar@{^{(}->}[ur] && 
M \ar@{^{(}->}[ur]
}
\end{displaymath}
Assume now that there is an extra occurrence of $L(e)$ in
$\mathcal{R}^1\mathrm{G}_{y}\Delta(e)$. This occurrence gives
us a homomorphism from $\Delta(e)$ to $\mathcal{R}^1\mathrm{G}_{y}\Delta(e)$,
which induces a non-zero homomorphism from $\Delta(e)$ to $M$.
Since $M$ embeds into $M'$ and the diagram commutes, 
our homomorphism defines a
homomorphism from $\Delta(e)$ to
$\mathcal{R}^1\mathrm{G}_{w_0}\Delta(e)$, which induces a
non-zero homomorphism from $\Delta(e)$ to $M'$. On the
other hand we know that
$[\mathcal{R}^1\mathrm{G}_{w_0}\Delta(e):L(e)]=
\overline{l}(w_0)$ by Lemma~\ref{l2.14}. From the previous
paragraph we also know that $[Z':L(e)]=\overline{l}(w_0)$. 
This gives us a contradiction and completes the proof
for the ungraded case. As we have mentioned above, the graded
version follows easily just tracking the grading.
\end{proof}

\begin{remark}\label{remextnew}
{\rm
Combined with Theorem~\ref{t2.02}\eqref{t2.02.2}, 
Theorem~\ref{t2.12} gives information about some 
higher $\mathrm{Ext}$-spaces. 
}
\end{remark}

\section{Modules with linear resolutions}\label{s4}

The category $\mathscr{LC}(\mathtt{P})$ realizes the  category 
of graded modules over the Koszul dual of $\mathtt{A}$ (which is 
isomorphic to $\mathtt{A}$ by \cite[Theorem~18]{So}). Verma modules
over $\mathtt{A}$
have linear projective resolutions. These resolutions, in turn, are
costandard objects in the category $\mathscr{LC}(\mathtt{P})$. 
In other words, this means that costandard modules are Koszul
dual to standard modules (but not vice versa). Analogously, 
since the Ringel dual  of $\mathtt{A}$ is isomorphic to 
$\mathtt{A}$ as well by \cite[Corollary~2.3]{So2},
costandard modules are also Ringel dual to standard modules
(but not vice versa). 

The category $\mathscr{LC}(\mathtt{T})$ realizes the 
category of graded modules over the Ringel dual of the 
Koszul dual of $\mathtt{A}$ 
(which is isomorphic to $\mathtt{A}$ by above).
Since the algebra $\mathtt{A}$ is standard Koszul (see 
\cite[Section~3]{ADL}), standard $\mathtt{A}$-modules admit 
linear tilting coresolutions and costandard $\mathtt{A}$-modules 
admit  linear tilting resolutions, see \cite[Theorem~7]{MO2}. 
In an analogy to the previous paragraph, from this one obtains
that both standard and costandard modules are Koszul-Ringel
self-dual. From \cite[Theorem~9]{MO2} it also follows that
simple and tilting $\mathtt{A}$-modules are Koszul-Ringel
dual to each other (now in the symmetric way). A natural
question then is: {\em Which other classes of modules can be
represented by linear complexes of tilting modules?}
(Such modules then in some sense ``live'' in the category $\mathscr{LC}(\mathtt{T})$). In this section we present
several classes of such modules. In particular,
quite surprizingly it turns our 
that all shuffled Verma modules have the above property.
In what follows we will use the term {\em tilting linearizable
modules} for those modules, which are isomorphic 
to some linear complexes of tilting modules in 
$\mathcal{D}^b(\mathtt{A}\mathrm{-fgmod})$.

\subsection{Shuffled Verma modules}\label{s4.1}

To start with we have to define graded lifts of shuffled Verma modules.
Let $\mathrm{T}_w:\mathtt{A}\mathrm{-gmod}\to\mathtt{A}\mathrm{-gmod}$
be the graded lift of $\mathrm{T}_w$, see \cite[Appendix]{MO2} or
\cite[page~28]{FKS}. We define the graded lifts of shuffled Verma 
modules as follows:
\begin{displaymath}
\mathtt{\Delta}(x,y)=\mathrm{T}_x\mathtt{\Delta}(y).
\end{displaymath}

\begin{theorem}\label{t16}
For every $x,y\in W$ the module $\mathtt{\Delta}(x,y)$
is tilting linearizable.
\end{theorem}

\begin{remark}\label{r17}
{\rm
The motivation for this statement is a compilation of several
results. \cite[Theorem~9]{MO2} and \cite[Corollary~14]{MO2}
say that in the category $\mathscr{LC}(\mathtt{T})\cong
\mathtt{A}\mathrm{-gmod}$ (which is a kind of 
``Koszul-Ringel dual'' to $\mathtt{A}\mathrm{-gmod}$) 
standard and costandard
$\mathtt{A}$-modules remain standard and costandard respectively,
and simple and tilting modules interchange. According to
\cite{AL}, shuffled Verma modules can be equivalently described 
using twisting and shuffling functors, the latter being
Koszul dual to each other by \cite[6.5]{MOS}. 
So it becomes natural to ask 
whether the set of shuffled Verma modules might be ``Koszul-Ringel 
self-dual''. The proof of Theorem~\ref{t16}, presented below, shows 
that this is indeed the case. Observe that it is very easy to see
on examples that this class is neither ``Ringel self-dual'' nor
``Koszul self-dual'' in general.
}
\end{remark}

\begin{proof}
The idea of the proof of Theorem~\ref{t16} is to compile the results
mentioned in Remark~\ref{r17}. The problem is to extend the 
``Koszul duality'' of shuffling and twisting functors from
\cite[6.5]{MOS} to the ``Koszul-Ringel duality'' of these functors.
For this we will need some notation.

Let $\mathrm{K}:\mathcal{D}^b(\mathtt{A}\mathrm{-gmod})\to 
\mathcal{D}^b(\mathscr{LC}(\mathtt{P}))$ denote the Koszul duality functor 
from \cite[5.4]{MOS} (restricted to bounded complexes).
Essentially this functor is given by taking the inner 
$\mathrm{Hom}$-functor with a direct sum of all indecomposable
projective objects from $\mathcal{D}^b(\mathscr{LC}(\mathtt{P}))$. 

By \cite[Theorem~2.2]{AS} and \cite[Theorem~6.6]{So2} for the functor
$\mathrm{T}_{w_0}$ we have that 
$\mathrm{T}_{w_0}:\mathcal{D}^b(\mathscr{LC}(\mathtt{P}))\to
\mathcal{D}^b(\mathscr{LC}(\mathtt{T}))$ is an equivalence, which sends
indecomposable projective objects from $\mathscr{LC}(\mathtt{P})$
to the corresponding indecomposable projective objects  from 
$\mathscr{LC}(\mathtt{T})$. This allows us to define the 
Koszul-Ringel duality functor
$\overline{\mathrm{K}}:\mathcal{D}^b(\mathtt{A}\mathrm{-gmod})\to 
\mathcal{D}^b(\mathscr{LC}(\mathtt{T}))$ as follows:
$\overline{\mathrm{K}}=\mathcal{L}\mathrm{T}_{w_0}\,\mathrm{K}$.

By \cite[6.4]{MOS}, translation and Zuckerman functors on
$\mathtt{A}\mathrm{-gmod}$ and $\mathscr{LC}(\mathtt{P})$ respectively 
are Koszul dual to each other with respect to the
Koszul duality $\mathrm{K}$. Since $\mathcal{L}\mathrm{T}_{w_0}$
commutes with translation functors by \cite[Theorem~3.2]{AS},
it follows that translation and Zuckerman functors on
$\mathtt{A}\mathrm{-gmod}$ and $\mathscr{LC}(\mathtt{T})$ respectively 
are Koszul-Ringel dual to each other with respect to the
Koszul-Ringel duality $\overline{\mathrm{K}}$.
Now, repeating the arguments from the proof of
\cite[Theorem~39]{MOS} one shows that twisting and shuffling
functors on $\mathtt{A}\mathrm{-gmod}$ and $\mathscr{LC}(\mathtt{T})$ respectively are Koszul-Ringel dual to each other with respect to the
Koszul-Ringel duality $\overline{\mathrm{K}}$. This means that
for any $w\in W$ we have
\begin{equation}\label{eqt16.1}
\mathcal{L}\mathrm{T}_w\cong\overline{\mathrm{K}}^{-1}\,
\mathcal{L}\overline{\mathrm{C}}_{w^{-1}}\,\overline{\mathrm{K}},
\end{equation}
where $\overline{\mathrm{C}}_{w^{-1}}$ denotes the corresponding 
shuffling  functor (see \cite{Ir3} and \cite[5.1]{MS2}).

The rest is now easy. Verma modules in 
$\mathtt{A}\mathrm{-gmod}$ and $\mathscr{LC}(\mathtt{T})$ correspond
via $\overline{\mathrm{K}}$ by \cite[Theorem~9]{MO2}.
Verma modules are acyclic for twisting functors by 
\cite[Theorem~2.2]{AS} and for shuffling functors by 
\cite[Proposition~5.3]{MS2}. Hence from \eqref{eqt16.1} for
$x,y\in W$ we have
\begin{displaymath}
\mathtt{\Delta}(x,y)=\mathrm{T}_x\mathtt{\Delta}(y)=
\overline{\mathrm{K}}^{-1}\,
\overline{\mathrm{C}}_{w^{-1}}\,\overline{\mathrm{K}}
\mathtt{\Delta}(y).
\end{displaymath}
Now, since the functor $\overline{\mathrm{C}}_{w^{-1}}$ is defined
already on $\mathscr{LC}(\mathtt{T})$, it follows that its 
value on $\overline{\mathrm{K}}^{-1}
\mathtt{\Delta}(y)\in \mathscr{LC}(\mathtt{T})$ is again an object
from $\mathscr{LC}(\mathtt{T})$. The necessary claim follows.
\end{proof}

\subsection{Standard modules in $\mathcal{O}_0^{\mathfrak{p}}$}\label{s4.2}

Let now $\mathfrak{p}\supset \mathfrak{h}\oplus\mathfrak{n}_+$ be
a parabolic subalgebra and $W^{\mathfrak{p}}$ the corresponding
parabolic subgroup of $W$. Let $\mathcal{O}_0^{\mathfrak{p}}$
denote the full subcategory of $\mathcal{O}_0$, consisting of
$U(\mathfrak{p})$-locally finite modules. Then simple 
objects of $\mathcal{O}_0$ have the form $L(w)$, where $w$ is the shortest
representative in a coset from $W^{\mathfrak{p}}\setminus W$.
We will denote the set of such representatives by $W(\mathfrak{p})$.
Let $A^{\mathfrak{p}}$ denote the  quotient of $A$ such that
$\mathcal{O}_0^{\mathfrak{p}}$ is equivalent to the category
of $A^{\mathfrak{p}}$-modules. Then $A^{\mathfrak{p}}$ is
quasi-hereditary (\cite{RC}) and inherits a positive
grading $\mathtt{A}^{\mathfrak{p}}$ for $\mathtt{A}$, with
respect to which it is  standard Koszul (\cite{BGS,ADL}). 
To indicate object of  $A^{\mathfrak{p}}$ we will add the
superscript $\mathfrak{p}$ to the standard notation.
As $\mathtt{A}^{\mathfrak{p}}$ is standard Koszul, the standard modules
$\mathtt{\Delta}^{\mathfrak{p}}(w)$, $w\in W(\mathfrak{p})$, 
have linear projective resolutions over $\mathtt{A}^{\mathfrak{p}}$. 
They also have linear tilting coresolutions over 
$\mathtt{A}^{\mathfrak{p}}$. Surprizingly enough, these
properties are preserved if one makes the step from 
$\mathtt{A}^{\mathfrak{p}}$ to $\mathtt{A}$.

\begin{proposition}\label{p3.01}
Let $w\in W(\mathfrak{p})$. Then, considered as an 
$\mathtt{A}$-module, the module $\mathtt{\Delta}^{\mathfrak{p}}(w)$
has a linear projective resolution and is tilting linearizable. 
\end{proposition}

\begin{proof}
The module $\Delta^{\mathfrak{p}}(w)$ is obtained via parabolic
induction (from $\mathfrak{p}$ to $\mathfrak{g}$) from a
simple finite-dimensional $\mathfrak{p}$-module. This
simple finite-dimensional $\mathfrak{p}$-module
has a BGG-resolution (over the Levi factor of $\mathfrak{p}$), 
which is obviously linear. The
parabolic induction then maps this BGG-resolution to a
linear resolution of $\Delta^{\mathfrak{p}}(w)$ by
standard modules over $\mathtt{A}$. Each standard
$\mathtt{A}$-module has a linear projective resolution and
a linear tilting coresolution. These resolutions
can be glued in the standard way to obtain linear 
projective resolution of $\Delta^{\mathfrak{p}}(w)$ and a
linear complex of tilting modules isomorphic to
$\Delta^{\mathfrak{p}}(w)$ respectively.
\end{proof}

\begin{remark}\label{rem4.09}
{\rm
I do not see any immediate
connection between the linear projective resolutions of
$\Delta^{\mathfrak{p}}(w)$ as $\mathtt{A}^{\mathfrak{p}}$-
and $\mathtt{A}$-modules.
}
\end{remark}

\begin{remark}\label{rem4.02}
{\rm
Applying $\mathrm{T}_{w_0}$ to the Verma resolution of
$\Delta^{\mathfrak{p}}(e)$ constructed in the proof of 
Proposition~\ref{p3.01} one obtains that
$\mathcal{L}\mathrm{T}_{w_0}\Delta^{\mathfrak{p}}(e)\cong
L(w_0^{\mathfrak{p}}w_0)[l(w_0^{\mathfrak{p}})]$. 
This allows one to compute the images of the simple 
modules $L(w_0^{\mathfrak{p}}w_0)$ under the (derived)
Ringel duality functor $\mathrm{Hom}_{A}(T,{}_-)$.
It is not clear how to compute these images
for other $L(x)$. This question reduces to understanding
the homology of the tilting objects in
$\mathscr{LC}(\mathtt{P})$ or of the projective objects
in $\mathscr{LC}(\mathtt{T})$.
}
\end{remark}

\begin{remark}\label{rem4.08}
{\em
Dually, costandard modules in a regular parabolic block
admit a linear injective coresolution, when viewed as
modules in the regular block of $\mathcal{O}$. Moreover, they
are also tilting linearizable.
}
\end{remark}

\subsection{Projective modules in $\mathcal{O}_0^{\mathfrak{p}}$}\label{s4.3}

\begin{proposition}\label{t4.05}
Let $w\in W(\mathfrak{p})$. Then, considered as an 
$\mathtt{A}$-module, the module $\mathtt{P}^{\mathfrak{p}}(w)$
has a linear projective resolution.
\end{proposition}

\begin{proof}
The module $\mathtt{P}^{\mathfrak{p}}(w)$ is obtained from 
$\mathtt{P}(w)$ by applying the $\mathfrak{p}$-Zuckerman functor.
Analogously to \cite[6.4]{MOS} one shows that the 
$\mathfrak{p}$-Zuckerman functor is Koszul dual to the
translation functor through the $W^{\mathfrak{p}}$-wall.
The latter functor preserves $\mathscr{LC}(\mathtt{P})$. Hence,
translating the simple object $\mathtt{P}(w)$ of
$\mathscr{LC}(\mathtt{P})$ through the $W^{\mathfrak{p}}$-wall
we will get a linear complex of projective modules, which
has only one non-zero homology, namely the one in the position
$0$, which is, moreover, isomorphic to 
$\mathtt{P}^{\mathfrak{p}}(w)$. The statement is proved.
\end{proof}

\begin{remark}\label{rem4.018}
{\em
Dually, injective modules in a regular parabolic block
of $\mathcal{O}$
admit linear injective coresolutions when
viewed as modules in $\mathcal{O}$.
}
\end{remark}

\subsection{Tilting modules in $\mathcal{O}_0^{\mathfrak{p}}$}\label{s4.4}

\begin{proposition}\label{t4.06}
Let $w\in W(\mathfrak{p})$. Then, considered as an 
$\mathtt{A}$-module, the module $\mathtt{T}^{\mathfrak{p}}(w)$
is tilting linearizable.
\end{proposition}

\begin{proof}
Apply $\mathcal{L}\mathrm{T}_{w_0}$ to the linear projective
resolution of $\mathtt{P}^{\mathfrak{p}}(x)$, $x\in W(\mathfrak{p})$,
constructed in Proposition~\ref{t4.05}, and follow the arguments 
of \cite[Proposition~4.4]{MS}.
\end{proof}

\begin{remark}\label{rem4.15}
{\em
From Propositions~\ref{t4.05} and \ref{t4.06} it follows that
projective tilting modules in $\mathcal{O}_0^{\mathfrak{p}}$
both admit a linear projective resolution in $\mathcal{O}$ and
are tilting linearizable. However, one has to note that a module
in  $\mathcal{O}_0^{\mathfrak{p}}$, which is at the same time
projective and tilting, has in the general case different graded 
lifts as a projective and as a tilting module. 
}
\end{remark}

\subsection{Some other classes of modules}\label{s4.5}

There are some other classes of modules, which are known to
have linear projective resolutions (respectively, which are 
tilting linearizable). In \cite[Proposition~4.1]{Ma}
it is shown that modules, obtained by translating standard modules
in singular blocks out of the wall, admit linear projective
resolutions. It is not difficult to show that they are also
tilting linearizable.
In \cite[Theorem~8.1 and Corollary~8.1]{Ma} it is shown that
one more class of modules (the ``wrong-sided'' analogue of 
modules, obtained by translating standard modules in singular blocks 
out of the wall) admits both a linear projective resolution
and a linear tilting coresolution.

The algebra $\mathtt{A}$ is an $\mathtt{A}\text{-}\mathtt{A}$ bimodule
and thus can be considered as an object of the category
$\mathcal{O}_0$ for the Lie algebra $\mathfrak{g}\times
\mathfrak{g}$ (this realization was used, in particular, in \cite{Ba}).
The hereditary chain of the quasi-hereditary algebra $\mathtt{A}$
is, by definition, a bimodule Verma flag for $\mathtt{A}$. From 
the natural grading on $\mathtt{A}$ we get that the heads of all the
Vermas occurring in this flag are concentrated in degree $0$.
Hence, we can glue linear projective resolutions (or linear
tilting coresolutions) of these Verma modules in the standard
way to obtain a linear projective resolution
(resp. a linear tilting coresolutions) of the bimodule 
$\mathtt{A}$. As a corollary one immediately obtains a
formula for computing Hochschild cohomology of $\mathtt{A}$
with coefficients in semi-simple modules.

\vspace{1cm}

\noindent
{\bf Acknowledgments}
The research was partially supported by The Royal Swedish Academy of 
Sciences, The Swe\-dish Research Council and The Swe\-dish Foundation
for International Cooperation in Research and Higher Education
(STINT).  An essential part of this paper was written during the
visit of the author to MPIM at Bonn. The financial support and
hospitality of MPIM is gratefully acknowledged. I would like to 
thank Catharina Stroppel and Olexandr Khomenko for several  very 
stimulating discussions and comments to the preliminary draft of
the paper, and to Jian-yi Shi for the information 
about two-sided cells he provided.

\vspace{0.5cm}

\noindent 
Department of Mathematics, Uppsala University, SE-751 06, Uppsala, SWEDEN,
e-mail: {\small \tt mazor@math.uu.se}, 
web: http://www.math.uu.se/$\tilde{\hspace{1mm}}$mazor/
\vspace{0.5cm}


\begin{thebibliography}{999}
\bibitem[ADL]{ADL}
{\em I.~{\'A}goston, V.~Dlab, E.~Luk{\'a}cs}, Quasi-hereditary 
extension algebras.  Algebr. Represent. Theory  {\bf 6}  (2003),  
no. 1, 97--117. 
\bibitem[AL]{AL}
{\em H.~Andersen, N.~Lauritzen}, Twisted Verma modules.  Studies 
in memory of Issai Schur (Chevaleret/Rehovot, 2000),  1--26, Progr. 
Math., {\bf 210}, Birkh{\"a}user Boston, Boston, MA, 2003.
\bibitem[AS]{AS}
{\em H.~Andersen, C.~Stroppel}, Twisting functors on $\mathcal{O}$. 
Represent. Theory {\bf 7} (2003), 681-699.
\bibitem[Ar]{Ar}
{\em S.~Arkhipov}, Algebraic construction of contragradient 
quasi-Verma modules in positive characteristic. in:
Representation theory of algebraic groups and quantum groups,
Adv. Stud. Pure Math. {\bf 40}, Math. Soc. Japan, Tokyo, 2004, pp. 27-68.
\bibitem[Ba1]{Ba}
{\em E.~Backelin}, The Hom-spaces between projective functors.  
Represent. Theory  {\bf 5}  (2001), 267-283.  
\bibitem[Ba2]{Ba2} 
{\em E.~Backelin}, Koszul duality for parabolic and singular 
category $\mathcal{O}$.  Represent. Theory  {\bf 3}  (1999), 139-152.
\bibitem[BGS]{BGS}
{\em A.~Beilinson, V.~Ginzburg, W.~Soergel}, Koszul duality patterns 
in representation theory.  J. Amer. Math. Soc.  {\bf 9}  (1996),  
no. 2, 473-527.
\bibitem[BG]{BG}
{\em I.~Bernstein, S.~Gelfand}, Tensor products of finite- 
and infinite-dimensional representations of semisimple Lie algebras.
Compositio Math.  {\bf 41}  (1980), no. 2, 245-285.
\bibitem[BGG1]{BGG}
{\em I.~Bernstein, I.~Gelfand, S.~Gelfand}, A certain category of 
${\mathfrak g}$-modules. Funkcional. Anal. i Prilozen.  {\bf 10}  
(1976), no. 2, 1-8.
\bibitem[BGG2]{BGG2}
{\em I.~Bernstein, I.~Gelfand, S.~Gelfand}, 
Differential operators on the base affine space and a study of 
${\mathfrak g}$-modules. in: Lie Groups Represent., Proc. Summer Sch. 
Bolyai Janos Math. Soc., Budapest 1971 (1975), 21-64.
\bibitem[BB]{BB}
{\em A.~Bj{\"o}rner, F.~Brenti}, Combinatorics of Coxeter groups. 
Graduate Texts in Mathematics, {\bf 231}. Springer, New York, 2005.
\bibitem[Bo]{Bo} 
{\em B.~Boe}, A counterexample to the Gabber-Joseph conjecture. 
Contemp. Math. {\bf 139} (1992), 1-3.
\bibitem[Ca]{Ca}
{\em K.~Carlin}, Extensions of Verma modules.  Trans. Amer. Math. 
Soc.  {\bf 294}  (1986),  no. 1, 29-43.
\bibitem[CPS]{CPS}
{\em E.~Cline, B.~Parshall, L.~Scott}, 
Finite dimensional algebras and highest weight categories. 
J. Reine Angew. Math. {\bf 391} (1988), 85-99.
\bibitem[CI]{CI}
{\em D.~Collingwood, R.~Irving}, A decomposition theorem for certain 
self-dual modules in the category ${\mathcal O}$. Duke Math. J. 
{\bf 58} (1989), no. 1, 89-102.
\bibitem[Di]{Di}
{\em J.~Dixmier}, Enveloping algebras. Revised reprint of the 
1977 translation. Graduate Studies in Mathematics, {\bf 11}, American 
Mathematical Society, Providence, RI, 1996.
\bibitem[DM]{DM}
{\em Yu.~Drozd, V.~Mazorchuk}, Koszul duality for extension 
algebras of standard modules, math.RT/0411528
\bibitem[FKS]{FKS}
{\em I.~Frenkel, M.~Khovanov, C.~Stroppel},
A categorification of finite-dimensional irreducible 
representations of quantum $sl(2)$ and their tensor products,
math.QA/0511467
\bibitem[GJ]{GJ}
{\em O.~Gabber, A.~Joseph}, Towards the Kazhdan-Lusztig conjecture. 
Ann. Sci. {\'E}c. Norm. Sup{\'e}r., IV. S{\'e}r. {\bf 14} (1981), 261-302.
\bibitem[Ha]{Ha}
{\em D.~Happel}, Triangulated categories in the representation theory 
of finite dimensional algebras. London Mathematical Society Lecture 
Note Series, 119. Cambridge etc.: Cambridge University Press. 1988.
\bibitem[Ir1]{Ir}
{\em R.~Irving}, Projective modules in the category 
$\mathcal{O}\sb S$: Loewy series.  Trans. Amer. Math. Soc.  
{\bf 291}  (1985),  no. 2, 733--754.
\bibitem[Ir2]{Ir2}
{\em R.~Irving},  A filtered category $\mathcal{O}\sb S$ 
and applications. Mem. Amer. Math. Soc. {\bf 83} (1990).
\bibitem[Ir3]{Ir3}
{\em R.~Irving}, Shuffled Verma modules and principal series 
modules over complex semisimple Lie algebras. J. London Math. Soc. 
(2) {\bf 48} (1993), no. 2, 263--277.
\bibitem[IS]{IS}
{\em R.~Irving, B.~Shelton}, Loewy series and simple projective modules 
in the category ${\mathcal O}_S$.  Pacific J. Math.  {\bf 132}  (1988),  
no. 2, 319-342.
\bibitem[Jo]{Jo}
{\em A.~Joseph}, The Enright functor on the 
Bernstein-Gelfand-Gelfand category $\mathcal{O}$.
Invent. Math. {\bf 67} (1982), no. 3, 423-445.
\bibitem[KL]{KL}
{\em D.~Kazhdan, G.~Lusztig}, Representations of Coxeter groups 
and Hecke algebras.  Invent. Math.  {\bf 53}  (1979), no. 2, 165-184. 
\bibitem[KM]{KM}
{\em O.~Khomenko, V.~Mazorchuk}, On Arkhipov's and Enright's 
functors. Math. Z. {\bf 249} (2005), no. 2, 357-386.
\bibitem[Lu1]{Lu}
{\em G.~Lusztig}, Cells in affine Weyl groups. Algebraic groups and 
related topics (Kyoto/Nagoya, 1983), 255-287, Adv. Stud. Pure 
Math., {\bf 6}, North-Holland, Amsterdam, 1985. 
\bibitem[Lu2]{Lu2}
{\em G.~Lusztig}, Cells in affine Weyl groups. II.
J. Algebra {\bf 109} (1987), no. {2}, 536-548.
\bibitem[Ma]{Ma}
{\em V.~Mazorchuk}, Applications of the category of linear 
complexes of tilting modules associated with the category 
$\mathcal{O}$, math.RT/0501220
\bibitem[MO1]{MO}
{\em V.~Mazorchuk, S.~Ovsienko}, Finitistic dimension of properly 
stratified algebras. Adv. Math. {\bf 186} (2004), no. 1, 251-265.
\bibitem[MO2]{MO2}
{\em V.~Mazorchuk, S.~Ovsienko}, A pairing in homology and the 
category of linear tilting complexes for a quasi-hereditary algebra,  
with an appendix by {\em C.~Stroppel}, J. Math. Kyoto Univ. 
{\bf 45} (2005), 711-741.
\bibitem[MOS]{MOS}
{\em V.~Mazorchuk, S.~Ovsienko, C.~Stroppel},
Quadratic duals, Koszul dual functors, and applications,
math.RT/0603475.
\bibitem[MS1]{MS}
{\em V.~Mazorchuk, C.~Stroppel}, Projective-injective modules, 
Serre functors and symmetric algebras, math.RT/0508119.
\bibitem[MS2]{MS2}
{\em V.~Mazorchuk, C.~Stroppel}, Translation and shuffling of 
projectively presentable modules and a categorification of a parabolic 
Hecke module.  Trans. Amer. Math. Soc.  {\bf 357}  (2005),  
no. 7, 2939-2973.
\bibitem[Ri]{Ri}
{\em C.~Ringel}, The category of modules with good filtrations 
over a quasi-hereditary algebra has almost split sequences. 
Math. Z. {\bf 208} (1991), No.2, 209-224.
\bibitem[RC]{RC}
{\em A.~Rocha-Caridi},  Splitting criteria for 
${\mathfrak g}$-modules induced from a parabolic and the 
Bernstein-Gelfand-Gelfand resolution of a finite-dimensional, 
irreducible ${\mathfrak g}$-module.  
Trans. Amer. Math. Soc.  {\bf 262}  (1980), no. 2, 335-366.
\bibitem[Sa]{Sa}
{\em B.~Sagan}, The symmetric group. Representations, combinatorial 
algorithms, and symmetric functions. Second edition. Graduate Texts 
in Mathematics, {\bf 203}. Springer-Verlag, New York, 2001.
\bibitem[So1]{So}
{\em W.~Soergel}, Kategorie $\mathcal{O}$, perverse Garben und 
Moduln {\"u}ber den Koinvarianten zur Weylgruppe. 
J. Amer. Math. Soc. {\bf 3} (1990), no. 2, 421-445.
\bibitem[So2]{So2}
{\em W.~Soergel}, Character formulas for tilting modules over 
Kac-Moody algebras.  Represent. Theory  {\bf 2}  (1998), 432-448.
\bibitem[St1]{St}
{\em C.~Stroppel},  Category $\mathcal{O}$: gradings and translation 
functors, J. Algebra {\bf 268}, no. 1, 301-326. 
\bibitem[St2]{St2}
{\em C.~Stroppel},  Homomorphisms and extensions of principal series
representations.  J. Lie Theory  {\bf 13}  (2003),  no. 1, 193-212.
\bibitem[Ve]{Ve}
{\em D.~Verma}, M{\"o}bius inversion for the Bruhat ordering on a 
Weyl group.  Ann. Sci. {\'E}cole Norm. Sup. (4)  {\bf 4}  (1971), 393-398.
\bibitem[Yu]{Yu}
{\em S.~Yuzvinsky}, Linear representations of posets, their cohomology 
and a bilinear form.  European J. Combin. {\bf 2}  (1981), no. 4, 385-397.
\bibitem[Zh]{Zh}
{\em B.~Zhu}, On characteristic modules of graded quasi-hereditary 
algebras.  Comm. Algebra  {\bf 32}  (2004),  no. 8, 2919-2928.
\end{thebibliography}
\end{document}